\documentclass[10pt]{article}

\usepackage{amsmath}

\usepackage{amssymb}

\usepackage{graphicx}
\def \s { \mathbf{\big|}}
\def \UN{\hbox{ \rm 1\hskip -3.2pt I}}
\def\R{\R}

\def\P{\mathbb{P}}
\def\E{\mathbb{E}} %

\def\R{\mathbb{R}}

\def \rhs {right hand side~}

\def\1{\mbox{I\hspace{-.6em}1}} 
\def\Var{\mathop{\rm Var}}

\newcommand{\findem}{\hfill\hbox{\hskip 4pt
\vrule width 5pt height 6pt depth 1.5pt}\vspace{.5cm}\par}

  \newtheorem{prop}{Proposition}
\newtheorem{cor}{Corollary}
\newtheorem{theo}{Theorem}
 \newtheorem{lem}{Lemma}

\makeindex
\begin{document}
\title{Some applications of Rice formulas to waves }

\author{ Jean-Marc Aza\"{\i}s \thanks { Universit\'e de Toulouse, IMT, LSP,
  F31062
Toulouse Cedex 9, France. Email: azais@cict.fr  } \and Jos\'e R. Le\'on
\thanks {Escuela de Matem\'atica. Facultad de Ciencias.
Universidad Central de Venezuela. A.P. 47197, Los Chaguaramos,
Caracas 1041-A, Venezuela. Email: jose.leon@ciens.ucv.ve} \and Mario
Wschebor
\thanks {Centro de Matem\'{a}tica. Facultad de
Ciencias. Universidad de la Rep\'{u}blica. Calle Igu\'{a} 4225.
11400. Montevideo. Uruguay. wschebor@cmat.edu.uy}}
\maketitle

\begin{abstract} We use Rice's formulas in order to compute the
moments of some level functionals which are linked to problems in
oceanography and optics. For instance, we consider the number of
specular points in one or two dimensions, the number of twinkles,
the distribution of normal angle of level curves and the number or
the length of dislocations in random wavefronts. We compute
expectations and in some cases, also second moments of such
functionals. Moments of order greater than one are more involved,
but one needs them whenever one wants to perform statistical
inference on some parameters in the model or to test the model
itself. In some cases we are able to use these computations to obtain
a Central Limit Theorem.
\end{abstract}


\section{Introduction}
Many problems in applied mathematics require to estimate the number
of points, the length, the volume and so on, of the level sets of a
random function $W(\mathbf{x})$, where $\mathbf{x}\in\R^d$, so that
one needs to compute the value of certain functionals of the
probability distribution of the size of the random set
$$
\mathcal{C}^W_A (\mathbf{u},\omega) := \{\mathbf{x}\in A:
W(\mathbf{x},\omega)=\mathbf{u}\},
$$
for some given $\mathbf{u}$.\\

Let us mention some examples which illustrate this general
situation:
\begin{itemize}
\item The number of times that a random process $\{X(t):t\in\R \}$ crosses the level
$u$:
$$
N_A^{X}(u)=\# \{ s \in A : X(s)=u\}.
$$
Generally speaking, the probability distribution of the random
variable $ N_A^{X}(u) $ is unknown, even for the simplest models of
the underlying process. However, there exist some formulas to
compute $\E( N_A^{X})$ and also higher order moments.

\item A particular case is the number of specular
points of a random curve or  a  random surface.\\
Consider first the case of a random curve. We take cartesian
coordinates $Oxz$ in the plane. A light source placed at $(0,h_1)$
emits a ray that is reflected at the point $(x,W(x))$ of the curve
and the reflected ray is registered by an observer placed at
$(0,h_2)$.

Using the equality between the angles of incidence and reflexion
with respect to the normal vector to the curve - i.e.
$N(x)=(-W'(x),1)$ - an elementary computation gives:
\begin{equation}\label{sp1}
  W'(x) = \frac{\alpha_2 r_1 - \alpha_1 r_2}{x(r_2-r_1)}
\end{equation}
where $ \alpha_i := h_i -W(x)$ and $r_i := \sqrt{ x^2 +\alpha_i^2 }
$,~~ i=1,2.\\

The points $(x,W(x))$ of the curve such that $x$ is a solution of
(\ref{sp1}) are called ``specular points''. We denote by $SP_1(A)$
the number of specular points such that $x\in A$, for each Borel
subset $A$ of the real line. One of our aims in this paper is to
study the probability distribution of $SP_1(A)$.

\item The following approximation, which turns out to be very accurate
in practice for ocean waves, was introduced long ago by
Longuet-Higgins (see [13] and
[14]):\\

Suppose that $h_1$ and $h_2$ are big with respect to $ W(x)$ and
$x$, then $ r_i = \alpha_i + x^2/( 2\alpha_i) + O( h_i^{-3} )$.
Then, (\ref{sp1}) can be approximated by
\begin{equation}\label{sp2}
  W'(x)   \simeq \frac{ x}{2 } \frac{\alpha_1 + \alpha_2}{\alpha_1 \alpha_2}
  \simeq \frac{ x}{2 }\frac{h_1 + h_2}{h_1 h_2} = k x,
\end{equation}
where
$$
k:=\frac{ 1}{2} \Big( \frac{1}{h_1} +\frac{1}{h_2} \Big).
$$

Denote $Y(x) := W'(x) -kx$ and $SP_2(A)$ the number of roots of
$Y(x)$ belonging to the set $A$, an approximation of $SP_1(A)$ under
this asymptotic. The first part of Section \ref{secspec} below will
be devoted to obtain some results on the distribution of the random
variable $SP_2(\R)$.

\item Consider now the same problem as above, but adding a time variable $t$, that is, $W$
becomes a random function parameterized by the pair $(x,t)$.
We denote $W_x,W_t,W_{xt},...$ the partial derivatives of $W$.\\

We use the Longuet-Higgins approximation (\ref{sp2}), so that the
approximate specular points at time $t$ are $(x,W(x,t))$ where
$$
W_x(x,t) =kx.
$$

Generally speaking, this equation defines a finite number of points
which move with time. The implicit function theorem, when it can be
applied, shows that the $x$-coordinate of a specular point moves  at
speed
$$
\frac{dx}{dt}=-\frac{W_{xt}}{W_{xx}-k}.
$$
The right-hand side diverges whenever $W_{xx}-k=0$, in which case a
flash appears and the point is called a ``twinkle''. We are
interested in the (random) number of flashes lying in a set $A$ of
space and in an interval $[0,T]$ of time. If we put:
\begin{equation}\label{f:twinkle}
\mathbf{Y}(x,t) := \left( \begin{array}{c}
  W_x(x,t) -kx \\
W_{xx}(x,t) -k
\end{array} \right).
\end{equation}
then, the number of twinkles  is:
$$
\mathcal{TW}(A,T)  := \sharp \{ (x,t) \in A \times [0,T]
:\mathbf{Y}(x,t) =0\}
$$

\item Let $W: Q \subset \R^d \to\R^{d'}$ with $d>d'$
be a random field and let us define the level set
$$
\mathcal{C}_Q^W(\mathbf{u})=\{\mathbf{x}\in Q:
W(\mathbf{x})=\mathbf{u}\}.
$$
Under certain general conditions this set is a $(d-d')$-dimensional
manifold but in any case, its $(d-d')$-dimensional Hausdorff measure
is well defined. We denote this measure by $\sigma_{d-d'}$. Our
interest will be to compute the mean of the $\sigma_{d-d'}$-measure
of this level set i.e.
$\E[\sigma_{d-d'}(\mathcal{C}_Q^W(\mathbf{u}))]$ as well as its
higher moments. It will be also of interest to compute:
$$
\E[\int_{\mathcal{C}_Q^W(\mathbf{u})}Y(s)d\sigma_{d-d'}(s)].
$$
where $Y(s)$ is some random field defined on the level set. Caba\~na
\cite{ca:ca}, Wschebor \cite{ws:ws} ($d'=1$) Aza\"{\i}s and Wschebor
\cite{aw}  and, in a weak form, Z\"ahle \cite{za:za} have studied
these types of formulas. See Theorems \ref{rice3} and \ref{rice4}.

\item Another interesting problem is the study of phase
singularities, dislocations of random wavefronts. They correspond to
lines of darkness, in light propagation, or threads of silence in
sound \cite{be1:be1}. In a mathematical framework they can be define
as the loci of points where the amplitude of waves vanishes. If we
represent the wave as
$$W(\mathbf{x},t)=\xi(\mathbf{x},t)+i\eta(\mathbf{x},t),\, \mbox{ where } \mathbf{x}\in \mathbb{R}^d$$ where
$\xi,\,\eta$ are independent homogenous Gaussian random fields the
dislocations are the intersection of the two random surfaces
$\xi(\mathbf{x},t)=0,\, \eta(\mathbf{x},t)=0$. We consider a fixed
time, for instance $t=0$. In the case $d=2$ we will study the
expectation of the following random variable
$$\#\{\mathbf{x}\in S:\,\xi(\mathbf{x},0)=\eta(\mathbf{x},0)=0\}.$$
In the case $d=3$ one important quantity is the length of the level
curve
$$\mathcal{L}\{\mathbf{x}\in
S:\,\xi(\mathbf{x},0)=\eta(\mathbf{x},0)=0\}.$$

\end{itemize}

All these situations are related to integral geometry. For a general
treatment of the basic theory, the classical reference is Federer's
``Geometric Measure Theory'' \cite{fe:fe}.

The aims of this paper are: 1) to re-formulate some known results in
a modern language or in the standard form of probability theory; 2)
to prove new results, such as computations in the exact models,
variance computations in cases in which only first moments have been
known, thus improving the statistical methods and 3) in some case,
obtain Central Limit Theorems.

The structure of the paper is the following: In Section
\ref{sectionrice} we review without proofs some formulas  for the
moments of the relevant random variables. In Section \ref{secspec}
we study expectation, variance and asymptotic behavior of specular
points. Section \ref{s:normal} is devoted to the study of the
distribution of the normal to the level curve.
 Section \ref{s:num} presents three numerical applications.
 Finally, in Section \ref{section4} we study dislocations of wavefronts
following  a paper by  Berry \& Dennis
\cite{be1:be1}.\\

\subsection*{Some additional notation and hypotheses}

$\lambda_d$ is Lebesgue measure in $\R^d$, $\sigma_{d'}(B)$ the
$d'$-dimensional Hausdorff measure of a Borel set $B$ and $M^T$ the
transpose of a matrix $M$. (const) is a positive constant whose
value may change from one occurrence to another.\\
If not otherwise stated, all random fields are assumed to be Gaussian and centered.

\section{Rice formulas}\label{sectionrice}

We give here a quick account of Rice formulas, which allow to
express the expectation and the higher moments of the size of level
sets of random fields by means of some integral formulas. The
simplest case occurs when both the dimension of the domain and the
range are equal to 1, for which the first results date back to Rice
\cite{rice} (see also Cram\'er and Leadbetter's book \cite{cr:le}).
When the dimension of the domain and the range are equal but bigger
than 1, the formula for the expectation is due to Adler \cite{a81}
for stationary random fields. For a general treatment of this
subject, the interested reader is referred to the book \cite{aw},
Chapters 3 and 6, where one can find proofs and details.
\begin{theo} [Expectation of the number of crossings, $d=d'=1$]\label{rice1}
Let $\mathcal{W}=\left\{ W(t) :t\in I \right\}, $ $I$ an interval in
the real line, be a Gaussian process having $\mathcal{C}^{1}$-paths.
Assume that $\Var(W(t))\neq 0$ for every $t\in I$. \\
Then:
\begin{equation} \label{riceg1}
\E\big( N_I^W(u)\big)  = \int_{I} \E\big( |W'(t)|\s W(t)=u \big)
p_{W(t)}(u) dt.
\end{equation}
\end{theo}

\begin{theo}[Higher moments of the number of crossings, $d=d'=1$]\label{ricek}
Let $m\geq 2$ be an integer. Assume that $\mathcal{W} $ satisfies
the hypotheses of Theorem \ref{rice1} and moreover, for any choice
of pairwise different parameter values $t_1,...,t_m \in I$ the joint
distribution of the $k$-random vector $(W(t_1),...,W(t_m))$ has a
density (which amounts to saying that its variance matrix is
non-singular). Then:
\begin{equation} \aligned \label{ricegk}
\E\big( N_I^W&(u)(N_I^W(u)-1)...(N_I^W(u)-m+1)\big)\\
&= \int_{I^m} \E\big( \prod_{j=1}^m|W'(t_j)|\s W(t_1)=...=W(t_m)=u
\big) p_{W(t_1),...,W(t_m)}(u,...,u) dt_1...dt_m.
\endaligned
\end{equation}
\end{theo}

Under certain conditions, the formulas in Theorems \ref{rice1} and
\ref{ricek} can be extended to non-Gaussian processes.

\bigskip

\begin{theo}[Expectation, $d=d'>1$]\label{Rice2}
Let $W:A\subset\R^d\to\R^d$ be a Gaussian random field, $A$ an open
set of  $\R^d$, $\mathbf{u}$ a fixed point in $\R^{d'}$ . Assume
that

\begin{itemize}
\item the sample paths  of $W$ are continuously
differentiable
\item for each $\mathbf{t}\in A$ the distribution of
$W(\mathbf{t})$ does not  degenerate
 \item $\mathbb{P}(\{\exists
\mathbf{t}\in A: W(\mathbf{t})=\mathbf{u}\,, \det(
W'(\mathbf{t}))=0\})=0$
\end{itemize}
Then for every Borel set  $B$ included in $A$
$$
\E\big(N_B^W(\mathbf{u})\big)=\int_B \E[|\det W'(\mathbf{t})|\s
W(\mathbf{t})=\mathbf{u}]p_{W(\mathbf{t})}(\mathbf{u})d\mathbf{t}.
$$
If $B$ is compact, both sides are finite.
\end{theo}
  The next
proposition provides sufficient conditions (which are mild) for the
third hypothesis in the above theorem to be verified (see again
\cite{aw}, Proposition 6.5).
\begin{prop}\label{prop1} Under the same conditions of the above theorem
one has
$$
\mathbb{P}(\{\exists \mathbf{t}\in A: W(\mathbf{t})=\mathbf{u}\,,
\det( W'(\mathbf{t}))=0\})=0
$$
if
\begin{itemize}
\item $p_{X(\mathbf{t})}(\mathbf{x})\le C$ for all $\mathbf{x}$ in
some neighborhood of $\mathbf{u}$,
\item at least one of the two
following conditions is satisfied\\
a) the trajectories of $W$ are twice continuously
differentiable\\
b)
$$
\alpha(\delta)=\sup_{x\in V(\mathbf{u})} \mathbb{P}\{|\det
W'(\mathbf{t})|<\delta \s W(\mathbf{t})=\mathbf{x}\}\to 0
$$
as $\delta\to 0$ where $V(\mathbf{u})$ is some neighborhood of
$\mathbf{u}$.
\end{itemize}
\end{prop}

\begin{theo}
[$m$-th factorial moment $d=d'>1$] \label{kricefield} Let $m\geq 2$
be an integer. Assume the same hypotheses as in
Theorem \ref{Rice2} except for (iii) that is replaced by\\

\noindent(iii') for $\mathbf{t}_{1},...,\mathbf{t}_{m}\in A$
distinct values of the parameter, the distribution of
$$
\big(W(\mathbf{t}_{1}),...,W(\mathbf{t}_{m})\big)
$$
does not degenerate in $(\R^{d})^{m}$.\\

Then for every Borel set B  contained in $A$,  one has
\begin{multline}\label{r13}
\E\left[ \big(N_B^W(\mathbf{u})\big)\big(N_B^W(\mathbf{u})-1\big)...\big(%
N_B^W(\mathbf{u})-m+1\big)\right]
 \\
=\int_{B^{m}}\E\Big( \prod_{j=1}^{m}|\det\big(W^{\prime
}(\mathbf{t}_{j})\big)| \s
W(\mathbf{t}_{1})=...=W(\mathbf{t}_{m})=u\Big)
\\
 p_{W(\mathbf{t}_{1}),...,W(\mathbf{t}_{m})}(\mathbf{u},...,\mathbf{u})d\mathbf{t}_{1}...d\mathbf{t}_{m},
\end{multline}
 where both sides may be infinite.
\end{theo}

\bigskip

When $d >d'$ we have the following formula :
\begin{theo}[Expectation of the geometric measure
of the level set. $d>d'$]\label{rice3}

Let $W : A\to \R^{d'}$ be a  Gaussian random field, $A$ an open
subset of $\R^d$, $d >d'$  and $\mathbf{u} \in \R^{d'}$ a fixed
point. Assume that:
\begin{itemize}
\item  Almost surely  the function  $\mathbf{t}\rightsquigarrow W(\mathbf{t})$ is of
class $\mathcal{C}^{1}$.
\item   For each $\mathbf{t}\in A$, $W(\mathbf{t})$ has a non-degenerate
distribution.
\item   $\P\{\exists \mathbf{t} \in A,W(\mathbf{t})=\mathbf{u},W'(\mathbf{t})~~\text{does not have
full rank} \}=0$
\end{itemize}
Then, for every Borel set B  contained in $A$,  one has
\begin{equation} \label{forlevel1}
\E \left( \sigma_{d-d'}(W,B)\right) =\int_{B}\E\left( \big[\det
\big(W'(\mathbf{t})(W'(\mathbf{t}))^T \big) \big]^{1/2} \s
W(\mathbf{t})=\mathbf{u}\right)
~p_{W(\mathbf{t})}(\mathbf{u})d\mathbf{t}.
\end{equation}
If B is compact, both sides in (\ref{forlevel1}) are finite.
\end{theo}

 The same kind of result holds true for integrals over the level set, as stated in the next
 theorem.

\begin{theo} [Expected integral on the level set]\label{rice4}

 Let $ W$ be a random field that verifies
the hypotheses of Theorem \ref{rice3}. Assume that for each
$\mathbf{t}\in A $ one has another random field $Y^\mathbf{t} : V
\to \R^{n}$, where $V$ is some topological space, verifying the
following conditions:
\begin{itemize}
  \item  $Y^\mathbf{t}(v)$ is a measurable function of $( \omega,\mathbf{t},v)$
  and almost surely, $(\mathbf{t},v)\rightsquigarrow Y^\mathbf{t}(v)$ is
  continuous.

  \item  For each $\mathbf{t}\in A$  the random process $(\mathbf{s},v) \to \big(
  W(\mathbf{s}), Y^\mathbf{t}(v) \big)$  defined on $W \times V$  is Gaussian.
\end{itemize}
Moreover, assume that $g: A\times \mathcal{C}(V, \R^{n}) \to \R $ is a
bounded function, which is continuous when one puts on
$\mathcal{C}(V, \R^{n})$ the topology of uniform convergence on
compact sets. Then, for each compact subset $B$ of $A$, one has
\begin{multline}\label{flevelpoids}
   \E\Big( \int_{B\cap W^{-1}(\mathbf{u})} g(t,Y^\mathbf{t}) \sigma_{d-d'}(W,d\mathbf{t}) \Big)
   \\
   =  \int_B \E \big( [\det(W'(t)(W'(t))^T) ]^{1/2} g(\mathbf{t},Y^\mathbf{t})\s Z(\mathbf{t})
   =\mathbf{u} \big) . p_{Z(\mathbf{t})} (\mathbf{u})
   d\mathbf{t}.
\end{multline}
\end{theo}

\section{Specular points and twinkles}\label{secspec}

\subsection{Number of roots} Let $W(\mathbf{t}):\R^d\to \R^d$ be a
zero mean stationary Gaussian field. If $W$ satisfies  the
conditions of Theorem \ref{Rice2} one has:
$$
\E\big(N_A^W(\mathbf{u})\big)= |A|\E[|\det
(W'(\mathbf{0}))|]p_{W(\mathbf{0})}(\mathbf{u}).
$$
where $ |A|$ denotes the Lebesgue measure of $A$.

 For $d=1$,
$N_A^W(u)$ is the number of crossings of the level $u$ and the
formula becomes
\begin{eqnarray}\label{Rice3}
\E\big(N_{[0,T]}^W(u)\big)=\frac
T\pi\sqrt{\frac{\lambda_2}{\lambda_0}}\,e^{-\frac{u^2}{2\lambda_0}},
\end{eqnarray}
 where
$$
\lambda_i= \int_0^\infty \lambda^i
d\mu(\lambda)\qquad i=0,2,4,\ldots,
$$
$\mu$ being  the spectral mesure
of $W$.

Formula (\ref{Rice3}) is in fact the one S.O. Rice wrote in the 40's see \cite{rice}.\\

\subsection{ Number of specular points}

We consider first the one-dimensional static case with the
longuet-Higgins approximation (\ref{sp2}) for the number of
specular points, that is:

$$
SP_2(I) = \# \{x \in I :  Y(x) = W'(x) -kx =0 \}
$$

 We assume that the Gaussian process $\{W(x):x\in \R \}$ has $
\mathcal{C}^2 $ paths and $\Var(W'(x))$ is constant equal to, say,
$v^2>0$. (This condition can always be obtained by means of an
appropriate non-random time change, the ``unit speed
transformation'') . Then Theorem 1 applies and
\begin{multline} \label{f:spec2}
\E (SP_2(I)) = \int_I\E(|Y'(x)| \s Y(x) =0) p_{Y(x)} (0) dx = \int_I
\E(|Y'(x)|)\frac{ 1}{v} \varphi(\frac{kx}{v}) dx
\\
= \int_I
G(-k,\sigma(x))\frac{1}{v} \varphi(\frac{kx}{v}) dv,
\end{multline}
 where $ \sigma^2(x)$ is the variance of $W''(x)$  and
 $ G(\mu, \sigma) :=  \E (|Z|)$,  $Z $ with  distribution $
 N(\mu,\sigma^2)$.\\

 For the second equality in (\ref{f:spec2}), in which we have erased the condition
 in the conditional expectation, take into account that
 since $\Var(W'(x))$ is constant, for each $x$ the random variables
 $W'(x) $ and $W''(x)$ are independent (differentiate under the
 expectation sign and use the basic properties of the Gaussian
 distribution).\\

 An elementary computation gives:
 $$
 G(\mu, \sigma) =\mu [2\Phi(\mu/\sigma )-1]+2 \sigma
 \varphi(\mu/\sigma),
$$
where  $ \varphi(.)$ and $\Phi(.)$  are respectively  the density
 and the cumulative distribution functions  of the standard Gaussian
 distribution. \medskip

 When the process $W(x)$ is also stationary, $ v^2 = \lambda_2$ and
 $ \sigma^2(x)$ is constant  equal to $ \lambda_4$. If we look at the  total number of
 specular points over the whole line, we get
\begin{equation} \label{f:spectot}
 \E(SP_2( \R)) =\frac{ G( k, \sqrt{\lambda_4})}{k}
\end{equation}
which is the result given by \cite{lo1:lo1} (part II, formula (2.14) page 846).
Note that this quantity
 is an increasing  function of $\frac{\sqrt{\lambda_4} )}{k}$.

  Since in the longuet-Higgins
approximation $k\approx 0$, one can write a Taylor expansion
having the form:
\begin{equation}\label{tayloresp}
\E(SP_2( \R)) \simeq
\sqrt{\frac{2\lambda_4}{\pi}}\frac{1}{k}\Big(
1+\frac{1}{2}\frac{k^2}{\lambda_4}+
\frac{1}{24}\frac{k^4}{\lambda_4^2}+...\Big)
\end{equation}

Let us turn to \textbf{the variance} of the number of specular
points, under some additional restrictions. First of all, we
assume for this computation that the given process $\{W(x):x\in \R
\}$ is stationary with covariance function\\
 $\E (W(x)W(y))=\Gamma
(x-y)$. $\Gamma$ is assumed to have enough regularity as to perform
the computations below, the precise requirements on it being given
in
the statement of Theorem \ref{asymptoticvar}.\\

Putting for short $S=SP_2( \R)$, we have:
\begin{equation}\label{vardesc}
\Var(S)=\E(S(S-1))+\E(S)-[\E(S)]^2
\end{equation}

The first term can be computed using Theorem \ref{ricek}:
\begin{align}
\E(S(S-1))&=\int\int_{\R^2}\E\big(|W''(x)-k||W''(y)-k|\s~
W'(x)=kx,W'(y)=ky
\big) \notag \\
&.p_{W'(x),W'(y)}(kx,ky)~dxdy\label{ricesp2}
\end{align}
where
\begin{equation}\label{densite}
p_{W'(x),W'(y)}(kx,ky)=\frac{1}{2\pi
\sqrt{\lambda_2^2-\Gamma''^2(x-y)}}\exp\Big[-\frac{1}{2}
\frac{k^2(\lambda_2x^2+2\Gamma''^2(x-y)xy+\lambda_2y^2)}{\lambda_2^2-\Gamma''^2(x-y)}
\Big],
 \end{equation}
  under the additional condition  that  the
density (\ref{densite}) does not degenerate for $x\neq y$.

 For the
conditional expectation in (\ref{ricesp2}) we perform a Gaussian
regression of $W''(x)$ (resp. $W''(y)$) on the pair $(W'(x),W'(y))$.
Putting $z=x-y$, we obtain:
\begin{equation}\label{regresionder2}
\aligned
&W''(x)=\theta_y(x)+a_y(x)W'(x)+b_y(x)W'(y)\\
&a_y(x)=-\frac{\Gamma'''(z)\Gamma''(z)}{\lambda_2^2-\Gamma''^2(z)}\\
&b_y(x)=-\frac{\lambda_2\Gamma'''(z)}{\lambda_2^2-\Gamma''^2(z)},
\endaligned
\end{equation}
where $\theta_y(x)$ is Gaussian centered, independent of
$(W'(x),W'(y))$. The regression of $W''(y) $ is obtained by
permuting $x$ and $y$.\\

The conditional expectation in (\ref{ricesp2}) can now be rewritten
as an unconditional expectation:
\begin{equation}\label{cond1}
 \E \Big\{ \Big|\theta_y(x)-k\Gamma'''(z)\big[
1+\frac{\Gamma''(z)x+\lambda_2y}{\lambda_2^2-\Gamma''^2(z)}\big]
\Big| \Big|\theta_x(y)-k\Gamma'''(-z)\big[
1+\frac{\Gamma''(-z)y+\lambda_2x}{\lambda_2^2-\Gamma''^2(z)}\big]
\Big|\Big\}
\end{equation}
Notice that the singularity on the diagonal $x=y$ is removable,
since a Taylor expansion shows that for $z\approx 0$:
\begin{equation}\label{eq1}
\Gamma^{\prime \prime \prime} (z)\big[
1+\frac{\Gamma''(z)x+\lambda_2y}{\lambda_2^2-\Gamma''^2(z)}\big]=
\frac{1}{2}\frac{\lambda_4}{\lambda_2}x \big(z+O(z^3) \big).
\end{equation}
One can check that
\begin{equation}\label{cond2}
\sigma ^2(z)=\E \big( (\theta_y(x))^2 \big)=\E \big(
(\theta_x(y))^2
\big)=\lambda_4-\frac{\lambda_2\Gamma'''^2(z)}{\lambda_2^2-\Gamma''^2(z)}
\end{equation}
and
\begin{equation}\label{cond3}
\E \big( \theta_y(x)\theta_x(y)
\big)=\Gamma^{(4)}(z)+\frac{\Gamma'''^2(z)\Gamma''(z)}{\lambda_2^2-\Gamma''^2(z)}.
\end{equation}

Moreover, if $\lambda _6<+\infty$, performing a Taylor expansion
one can show that as $z\approx 0$ one has
\begin{equation}\label{eq2}
\sigma ^2(z)\approx \frac{1}{4}\frac{\lambda _2\lambda _6-\lambda
_4^2}{\lambda _2}z^2
\end{equation}
and it follows that the singularity at the diagonal of the
integrand
in the right-hand side of (\ref{ricesp2}) is also removable.\\

We will make use of the following auxiliary statement that we state
as a lemma for further reference.The proof requires some
calculations, but is elementary and we skip it. The value of
$H(\rho;0,0)$ below can be found for example in \cite{cr:le}, p.
211-212.
\begin{lem}\label{esproabsv}
Let
$$
H(\rho;\mu,\nu )=\E (|\xi+\mu||\eta+\nu|)
$$
where the pair $(\xi,\eta)$ is centered Gaussian, $\E(\xi
^2)=\E(\eta ^2)=1,~~\E(\xi\eta)=\rho.$\\

Then,
$$
H(\rho;\mu,\nu)= H(\rho;0,0)+R_2(\rho;\mu,\nu)
$$
where
$$
H(\rho;0,0)=\frac{2}{\pi}\sqrt{1-\rho ^2}+\frac{2\rho}{\pi}\arctan
\frac{\rho}{\sqrt{1-\rho ^2}},
$$
and
$$
|R_2(\rho;\mu,\nu)|\leq 3 (\mu^2+\nu^2)
$$
if $\mu^2+\nu^2 \leq 1$ and $0\leq \rho \leq 1$.
\end{lem}

In the next theorem we compute the equivalent of the variance of the
number of specular points, under certain hypotheses on the random
process and with the longuet-Higgins asymptotic. This result is new
and useful for estimation purposes since it implies that, as
$k\rightarrow 0$, the coefficient of variation of the random
variable $S$ tends to zero at a known speed. Moreover, it will also
appear in a natural way when normalizing $S$ to obtain a Central
Limit Theorem.

\begin{theo}\label{asymptoticvar}
Assume that the centered Gaussian stationary process
$~\mathcal{W}=\{W(x):x\in \R \}$ is $\delta-$dependent, that is,
$\Gamma (z)=0$ if $|z|>\delta$, and that it has
$\mathcal{C}^4$-paths. Then, as $k\rightarrow 0$ we have:
\begin{equation}\label{eqvar}
\Var (S)=\theta \frac{1}{k}+O(1).
\end{equation}
where
$$
\theta=\Big(\frac{J}{\sqrt{2}}+\sqrt{\frac{2\lambda
_4}{\pi}}-\frac{2\delta \lambda _4}{\sqrt{\pi^3 \lambda _2}} \Big),
$$

\begin{equation}\label{intJ}
J=\int_{-\delta}^{+\delta}\frac{\sigma^2(z)H\big(\rho(z);0,0)
\big)}{\sqrt{2\pi(\lambda_2+\Gamma''(z))}}dz,
\end{equation}
the functions $H$ and $\sigma^2(z)$ have already been defined above,
and
$$
\rho(z)=\frac{1}{\sigma ^2(z)}\Big[ \Gamma ^{(4)}(z)+\frac{\Gamma
'''(z)^2\Gamma ''(z)}{\lambda _2^2-\Gamma''^2(z)}\Big].
$$
\end{theo}

\textbf{Remarks on the statement}.
\begin{itemize}
\item The assumption that the paths of the process are of class
$\mathcal{C}^4 $ imply that $\lambda _8<\infty$. This is well-known
for Gaussian stationary processes (see for example \cite{cr:le}).
\item Notice that since the process is $\delta$-dependent, it is also
$\delta'$-dependent for any $\delta'>\delta$. It is easy to verify
that when computing with such a $\delta'$ instead of $\delta$ one
gets the same value for $\theta$.
\item One can replace the $\delta$-dependence
by some weaker mixing condition, such as
$$
\big|\Gamma ^{(i)}(z) \big|~ \leq ~ (const)(1+|z|)^{-\alpha}~~~(0\leq
i \leq 4)
$$
for some $\alpha >1$, in which case the value of $\theta$ should be
replaced by:
$$
\theta=\sqrt{\frac{2\lambda_4}{\pi}}+\frac{1}{\sqrt{\pi}}\int_{-\infty}^{+\infty}
\Big[ \frac{\sigma ^2(z)H(\rho (z);0,0)}{2\sqrt{\lambda_2+\Gamma
''(z)}}-\frac{1}{\pi}\frac{\lambda _4}{\sqrt{\lambda _2}}\Big]dz.
$$
The proof of this extension can be performed following the same
lines as the one we give below, with some additional computations.
\end{itemize}

\bigskip

 \emph{Proof of the Theorem}:
We use the notations and computations preceding the statement of
the
theorem.\\
Divide the integral on the right-hand side of (\ref{ricesp2}) into
two parts, according as $|x-y|>\delta$ or $|x-y|\leq \delta$, i.e.

\begin{equation}\label{dividiren2}
\E (S(S-1))=\int\int_{|x-y|>\delta}...+\int\int_{|x-y|\leq
\delta}...=I_1+I_2.
\end{equation}

In the first term, the $\delta-$dependence of the process implies
that one can factorize the conditional expectation and the density
in the integrand. Taking into account that for each $x\in \R$, the
random variables $W''(x)$ and $W'(x)$ are independent, we obtain for
$I_1$:
\begin{equation*}
I_1=\int\int_{|x-y|>\delta}\E \big(|W''(x)-k| \big)\E
\big(|W''(y)-k| \big)p_{W'(x)}(kx)p_{W'(y)}(ky)dxdy.
\end{equation*}
On the other hand, we know that $W'(x)$ (resp. $W''(x)$) is
centered normal with variance $\lambda_2$ (resp. $\lambda_4$).
Hence:
$$
I_1=\big[G(k,\sqrt{\lambda_4})
\big]^2\int\int_{|x-y|>\delta}\frac{1}{2\pi \lambda_2}\exp\Big[
-\frac{1}{2}\frac{k^2(x^2+y^2)}{\lambda_2} \Big]dxdy,
$$
To compute the integral on the right-hand side, notice that the
integral over the whole $x,y$ plane is equal to $1/k^2$ so that it
suffices to compute the integral over the set ${|x-y|\leq \delta}$.
Changing variables, this last  one is equal to

\begin{equation*}\aligned
&\int_{-\infty}^{+\infty}dx\int_{x-\delta}^{x+\delta}\frac{1}{2\pi
\lambda _2}\exp\Big[ -\frac{1}{2}\frac{k^2(x^2+y^2)}{\lambda_2}
\Big]dy\\
&=\frac{1}{2\pi k^2}\int_{-\infty}^{+\infty} e^{-\frac{1}{2}u^2}du
\int_{u-\frac{k\delta}{\sqrt{\lambda_2}}}^{u+\frac{k\delta}{\sqrt{\lambda_2}}}
e^{-\frac{1}{2}v^2}dv\\
&=\frac{\delta}{k\sqrt{\lambda_2\pi}}+O(1),
\endaligned
\end{equation*}
where the last term is bounded if $k$ is bounded (in fact, remember
that we are considering an approximation in which $k\approx 0$). So,
we can conclude that:
$$
\int\int_{|x-y|>\delta}\frac{1}{2\pi \lambda_2}\exp\Big[
-\frac{1}{2}\frac{k^2(x^2+y^2)}{\lambda_2}
\Big]dxdy=\frac{1}{k^2}-\frac{\delta}{k\sqrt{\lambda_2\pi}}+O(1)
$$
Replacing in the formula for $I_1$ and performing a Taylor
expansion, we get:
\begin{equation}\label{i1}
I_1=\frac{2\lambda_4}{\pi}\Big[
\frac{1}{k^2}-\frac{\delta}{k\sqrt{\lambda_2 \pi}}+O(1)\Big].
\end{equation}

Let us now turn to $I_2$.\\

Using Lemma \ref{esproabsv} and the equivalences (\ref{eq1}) and
(\ref{eq2}), whenever $|z|=|x-y|\leq \delta$, the integrand on the
right-hand side of (\ref{ricesp2}) is bounded by
$$
(const)\big[ H(\rho (z);0,0)+k^2(x^2+y^2)\big].
$$

We divide the integral $I_2$ into two parts:\\

First, on the set $\{(x,y):|x|\leq 2\delta, |x-y|\leq
\delta \}$ the integral is clearly bounded by some constant.\\

Second, we consider the integral on the set $\{(x,y):x> 2\delta,
|x-y|\leq \delta \}$. (The symmetric case, replacing $x> 2\delta$ by
$x<- 2\delta$ is
similar,that is the reason for the factor $2$ in what follows).\\

We have (recall that $z=x-y$):

\begin{equation*}\aligned
I_2=O(1)+2\int&\int_{|x-y|\leq \delta,x>2\delta}\sigma^2(z)\Big[
H\big(\rho (z);0,0 \big)+R_2\big(\rho(z);\mu,\nu \big)\Big]\\
&\times \frac{1}{2\pi \sqrt{\lambda _2^2-\Gamma ''^2(z)}}\exp
\Big[-\frac{1}{2}\frac{k^2(\lambda _2x^2+2\Gamma ''(x-y)xy+\lambda
_2y^2)}{\lambda _2^2-\Gamma''^2(x-y)}\Big]dxdy
\endaligned
\end{equation*}
which can be rewritten as:
\begin{equation*}\aligned
I_2=O(1)+2\int_{-\delta} ^{\delta}&\sigma^2(z)\Big[
H\big(\rho (z);0,0 \big)+R_2\big(\rho(z);\mu,\nu \big)\Big]\\
&\times \frac{1}{\sqrt{2\pi (\lambda _2+\Gamma ''(z))}}\exp
\Big[-\frac{1}{2}\frac{k^2z^2}{\lambda_2-\Gamma''(z)}
\big(\frac{\lambda_2}{\lambda_2+\Gamma''(z)}-\frac{1}{2}
\big)\Big]dz\\
&\times \int_{2\delta}^{+\infty}\frac{1}{\sqrt{2\pi (\lambda
_2-\Gamma ''(z))}}\exp\Big[-k^2\frac{(x-z/2)^2}{\lambda _2-\Gamma
''(z))} \Big]dx
\endaligned
\end{equation*}

In the inner integral we perform the change of variables
$$
\tau =\frac{\sqrt{2} k(x-z/2)}{\sqrt{\lambda _2-\Gamma ''(z)}}
$$
so that it becomes:
\begin{equation}\label{intint}
\frac{1}{k\sqrt{2}}\int_{\tau _0}^{+\infty}\frac{1}{\sqrt{2\pi
}}\exp\big( -\frac{1}{2}\tau ^2
\big)d\tau~=~\frac{1}{2\sqrt{2}}\frac{1}{k}+O(1)
\end{equation}
where
$\tau_0=2\sqrt{2}k(2\delta-z/2)/\sqrt{\lambda_2-\Gamma''(z)}$.\\

\noindent Notice that $O(1)$ in (\ref{intint}) is uniformly bounded,
independently of $k$ and $z$, since the hypotheses on the process
imply that $\lambda_2-\Gamma''(z)$ is bounded below by a positive
number, for
all $z$.\\

\noindent We can now replace in the expression for $I_2$ and we obtain
\begin{equation}\label{finali2}
I_2=O(1)+\frac{J}{k\sqrt{2}}.
\end{equation}
To finish, put together (\ref{finali2}) with (\ref{i1}),
(\ref{dividiren2}), (\ref{vardesc}) and (\ref{tayloresp}).  \findem

\begin{cor}\label{corcoefvar}
Under the conditions of Theorem \ref{asymptoticvar}, as
$k\rightarrow 0 $:
$$
\frac{\sqrt{\Var(S)}}{\E(S)}\approx \sqrt{\theta k}.
$$
\end{cor}

The proof follows immediately from the Theorem and the value
of the expectation.\\

 The computations made in this section are in
close relation with the two results of Theorem $4$ in Kratz and
Le\'on \cite{kra:kra}. In this paper the random variable $SP_2(I)$
is expanded in the Wiener-Hermite Chaos. The aforementioned
expansion yields the same formula for the expectation and allows
obtaining also a formula for the variance. However, this expansion
is difficult to manipulate in order to get the result of Theorem
\ref{asymptoticvar}.

\bigskip

 Let us now turn to the Central Limit Theorem.
 \begin{theo}\label{tcl}
Assume that the process $\mathcal{W} $ satisfies the hypotheses of
Theorem \ref{asymptoticvar}. In addition, we will assume that the
fourth moment of the number of approximate specular points on an
interval having length equal to $1$  is bounded uniformly in $k$,
that is
\begin{equation}\label{mom4}
\E \big(\big[SP_2([0,1]) \big]^4\big)\leq (const)
\end{equation}

Then, as $ k\rightarrow 0$,
\begin{equation*}
\frac {S - \sqrt{\frac{2\lambda _4}{\pi}}\frac{1}{k}}{\sqrt{\theta
/k}} \Rightarrow~~N(0,1),
\end{equation*}
where $\Rightarrow $ denotes convergence in distribution.
  \end{theo}

\textbf{Remark.}\\

One can give conditions for the added hypothesis (\ref{mom4}) to
hold true, which require some additional regularity for the process.
Even though they are not nice, they are not costly from the point of
view of physical models. For example, either one of the following
conditions imply (\ref{mom4}):

   \begin{itemize}
     \item The paths $x\rightsquigarrow W(x)$ are of class $\mathcal{C}^{11}$. (Use Theorem 3.6 of \cite{aw}
     with $m=4$, applied to the random process $\{W'(x):x\in \R \}$.
     See also \cite {nu:ws}).
     \item  The paths $x\rightsquigarrow W(x)$ are of class $\mathcal{C}^{9}$ and the support of the spectral
     measure has an accumulation  point: apply
     Exercice  3.4 of  \cite{aw} to get the non-degeneracy condition, Proposition  5.10 of
     \cite{aw} and Rice formula (Theorem \ref{ricek}) to get
     that the fourth moment of the number of zeros of $W''(x)$ is
     bounded.
   \end{itemize}

   {\em Proof of the Theorem. }
Let $\alpha$ and $\beta$ be real numbers satisfying the conditions
   $1/2<\alpha<1$, $\alpha + \beta >1$, $2\alpha + \beta <2$.
   It suffices to prove the convergence as $k$ takes values on a
   sequence of positive numbers tending to $0$. To keep in mind that the parameter is $k$,
   we use the notation
   $$
   S(k):=S=SP_2(\R)
   $$

   Choose $k$ small enough, so that $k^{-\alpha} >2$ and define
   the sets of disjoint intervals, for $j=0, \pm 1,\ldots, \pm [k
^{-\beta}]$:

\begin{align*}
  U_j^k &= \big((j-1) [k ^{-\alpha}]\delta+\delta /2, j  [k ^{-\alpha}]\delta -\delta /2\big),
\\
I_j^k&=\big[j[k ^{-\alpha}]\delta-\delta /2,j[k
^{-\alpha}]\delta+\delta /2 \big].
\end{align*}
$[.]$ denotes integer part.\\

Notice that each interval $U_j^k$ has length $[k ^{-\alpha}]\delta
-\delta$ and that two neighboring intervals $U_j^k$ are separated by
an interval of length $\delta$. So, the $\delta$-dependence of the
process implies that the random variables $SP_2(U_j^k),~j=0, \pm
1,\ldots, \pm [k^{-\beta}]$ are independent. A similar argument
applies to $SP_2(I_j^k),~j=0, \pm
1,\ldots, \pm [k^{-\beta}]$.\\

We denote:
 $$
   T(k) = \sum_{|j|\leq [k ^{-\beta}]} SP_2( U_j^k),
   $$
Denote
$$
V_k = \big(\Var( S(k)) \big) ^{-1/2} \approx \sqrt{k/\theta}
$$
where the equivalence is due to Theorem \ref{asymptoticvar}.\\

We give the proof in two steps, which easily imply the statement. In
the first one, we prove that
$$
 V_k [ S(k) -T(k)]
$$
 tends to $0$ in the $L^2$ of the underlying probability space.

In the second step we prove that
$$
 V_kT(k)
$$
is asymptotically standard normal.\\

\textbf{Step 1}. We prove first that $V_k [S(k) -T(k)] $ tends to
$0$ in $L^1$. Since it is non-negative, it suffices to show that its
expectation tends to zero. We have:
$$
S(k) -T(k)= \sum _{|j|< [k ^{-\beta}]}SP_2 (I_j^k)+Z_1+Z_2
$$
where\\

$Z_1=SP_2 \big(-\infty,- [k ^{-\beta}]. [k ^{-\alpha}]\delta +\delta
/2 \big)$,\\

$Z_2=SP_2 \big( [k ^{-\beta}]. [k ^{-\alpha}]\delta-\delta
/2,+\infty)\big).$\\

\noindent Using the  fact that $\E\big(SP^k_2(I)\big) \leq (const) \int_I
\varphi( kx/\sqrt{\lambda_2}) dx, $  one can show that
 $$ V_k \E(S(k) -T(k)) \leq  (const)k^{1/2} \Big[ \sum_{\ell  =0}^{+\infty} \varphi \big(\frac{\ell [k^{-\alpha}]k\delta }{ \sqrt{\lambda_2}}\big)
 + \int_{ [k^{-\alpha}][k^{-\beta}]\delta} ^{+\infty} \varphi( kx/\sqrt{\lambda_2}) dx \Big].
 $$
which tends to zero as a consequence of the choice of $\alpha$ and
$\beta$.\\

It suffices to prove that $V_k^2 \Var \big( S(k) -T(k) \big)
\rightarrow 0 $ as $k\rightarrow 0$. Using independence:

\begin{equation}\aligned \label{varsuma}
\Var \big( S(k) -T(k) \big)&= \sum _{|j|< [k ^{-\beta}]}\Var \big(
SP_2 (I_j^k) \big) + \Var(Z_1)+ \Var (Z_2 )\\
&\leq \sum _{|j|< [k ^{-\beta}]}\E \big( SP_2 (I_j^k)(SP_2
(I_j^k)-1) \big)\\
&+ \E(Z_1(Z_1-1))+ \E (Z_2(Z_2-1) )+\E \big( S(k)-T(k) \big).
\endaligned
\end{equation}

We already know that $V_k^2~\E \big( S(k)-T(k) \big)\rightarrow
0.$ Using the hypotheses of the theorem, since each $I_j^k$ can be covered by a fixed  number of intervals of size one, we know that
$\E \big( SP_2 (I_j^k)(SP_2 (I_j^k)-1) \big) $ is
bounded by a constant which does not depend on $k$ and $j$.
We can write
$$
 V_k^2\sum _{|j|< [k ^{-\beta}]}\E \big( SP_2
(I_j^k)(SP_2 (I_j^k)-1) \big)\leq (const)k^{1-\beta}
$$
which tends to zero because of the choice of $\beta$.
The remaining two terms  can be bounded by calculations similar to those of the proof of Theorem \ref{asymptoticvar}.

\medskip

 {\bf Step 2.} $~~~T(k)$ is a sum of independent
 but not equi-distributed random variables. To prove it satisfies a Central Limit Theorem,
 we use a Lyapunov condition based of fourth moments. Set:
$$
M_j^m := \E \Big\{\big[ SP_2(U_j^k) -\E\big( SP_2(U_j^k)\big)
\big]^m\Big\}
$$
For the Lyapunov condition it suffices to verify that
\begin{equation}\label{lyap}
\Sigma^{-4}\sum_{|j|\leq [k ^{-\beta}]} M_j^4  \to 0 ~~~\text{as} ~~
k\rightarrow 0,
\end{equation}
where
$$
\Sigma^2 := \sum_{|j|\leq [k ^{-\beta}]} M_j^2.
$$
To prove (\ref{lyap}), let us  partition each interval $U_j^k$ into
$ p=[k^{-\alpha}]-1$ intervals $I_1,...I_p$ of equal size $\delta$.
We have
\begin{equation}\label{sumamom4}
\E\big(SP_1+\cdots+SP_p)^4 = \sum_{1\leq i_1,i_2,i_3,i_4 \leq p} \E
\big( SP_{i_1}SP_{i_2} SP_{i_3} SP_{i_4}\big),
\end{equation}
where $ SP_i$ stands for $SP_2(I_i) -\E\big( SP_2(I_i)\big)$ Since
the size of all the intervals
 is equal to $\delta$ and given the finiteness of fourth moments in the hypothesis, it follows that
 $\E \big( SP_{i_1}SP_{i_2} SP_{i_3} SP_{i_4}\big)$ is bounded.\\

On the other hand, notice that the number of terms which do not
vanish in the sum of the right-hand side of (\ref{sumamom4}) is
$\mathcal{ O}(p^2 )$. In fact, if one of the indices in
$(i_1,i_2,i_3,i_4)$ differs more than $1$ from all the other, then
$\E \big( SP_{i_1}SP_{i_2} SP_{i_3} SP_{i_4}\big)$ vanishes. Hence,
$$
\E \Big[ SP_2(U_j^k) -\E\big( SP_2(U_j^k)\big) \Big]^4 \leq (const)
k^{-2\alpha}
$$
so that $\sum_{|j|\leq [k^{-\beta}]} M_j^4 = \mathcal{O}
(k^{-2\alpha} k^{-\beta}).$ The inequality $2 \alpha +\beta <2 $
implies Lyapunov condition. \findem

 \subsection{Number of specular points without approximation}

 We turn now to the computation of the expectation of the number
  of specular points  $SP_1(I)$ defined by
 (\ref{sp1}).
 This number of specular points is equal to the number of zeros of
 the process
 $$
 Z (x) := W'(x) -m_1(x,W(x)) =0,
$$
where
$$
m_1(x,w)=\frac{x^2-(h_1-w)(h_2-w)+\sqrt{[x^2+(h_1-w)^2][x^2+(h_2-w)^2]}}{x(h_1+h_2-2w)}.
$$
Assume that the process $\{W(x):x\in \R \}$ is Gaussian, centered,
stationary, with $\lambda_0=1$. The process $Z(t)$  is not
Gaussian and we must use a generalization of Theorem \ref{rice1},
namely Theorem 3.2 of \cite{aw} to get
\begin{multline}\label{spec1}
 \E \big( SP_1([a,b])\big)=\int_a^b dx \int _{-\infty}
^{+\infty}\E \big( |Z'(x)|\s Z(x)=0,W(x)=w \big)\\
.\frac{1}{\sqrt{2 \pi}} e^{-\frac{w^2}{2}} \frac {1}{\sqrt{2\pi
\lambda _2}}e^{-\frac {m_1^2(x,w)}{2 \lambda _2}} dw.
\end{multline}
For the conditional expectation in (\ref{spec1}), notice that
$$
Z'(x)=W''(x)-\frac{\partial m_1}{\partial
x}(x,W(x))-\frac{\partial m_1}{\partial w}(x,W(x))W'(x),
$$
so that under the condition,
$$
Z'(x)=W''(x)-K(x,w),~~ \text{where} ~~K(x,w)=\frac{\partial
m_1}{\partial x}(x,w))+\frac{\partial m_1}{\partial
w}(x,w))m_1(x,w).
$$

Using that for each $x$,  $W''(x)$ and $W'(x)$ are independent
random variables and performing a Gaussian regression of $W''(x)$
on $W(x)$, we can write (\ref{spec1}) in the form:

\begin{multline}\label{spec2}
\E \big( SP_1([a,b])\big)\\
=\int_a^b dx \int _{-\infty} ^{+\infty}\E
\big( |\zeta-\lambda _2w -K(x,w)| \big) \frac{1}{2 \pi
\sqrt{\lambda_2}} \exp \Big(-\frac{1}{2}(w^2+\frac {m_1^2(x,w)}{\lambda
_2}) \Big) dw.
\end{multline}
where $\zeta$ is centered Gaussian with variance
$\lambda_4-\lambda_2^2$.
Formula (\ref{spec2}) can still be rewritten as:
\begin{multline}\label{spec3}
\E \big( SP_1([a,b])\big)\\
=\frac{1}{2
\pi}\sqrt{\frac{\lambda_4-\lambda_2^2}{\lambda_2}}\int_a^b dx \int
_{-\infty} ^{+\infty}
G(m,1)
 \exp\Big(-\frac{1}{2}(w^2+\frac {m_1^2(x,w)}{\lambda
_2})\Big) dw,
\end{multline}
where
$$
m=m(x,w)=\frac{\lambda_2w+K(x,w)}{\sqrt{\lambda_4-\lambda_2^2}}.
$$

Notice that in (\ref{spec3}), the integral is convergent as
$a\rightarrow -\infty,~b\rightarrow +\infty$ and that this formula
is well-adapted to numerical approximation.\\

\subsection{Number of twinkles}

We give a proof of a result stated in
\cite{lo1:lo1} (part III pages 852-853).

We consider $\mathbf{ Y}(x,t)$ defined by (\ref{f:twinkle}) and we
limit ourselves to the case in which $ W(x,t)$ is centered  and
stationary. If $\mathbf{Y}$ satisfies the conditions of Theorem
\ref{Rice2}, by stationarity we get
\begin{equation} \label{f:twin}
\E\big( \mathcal{TW}(I,T)\big) = T \int_I \E\big( |\det
\mathbf{Y}'(x,t) | \s \mathbf{Y}(x,t) =0 \big)
p_{\mathbf{Y}(x,t)}(0)dx.
\end{equation}
Since $W_{xx}$ and $W_x$ are independent  with respective variances
$$
\lambda_{40} = \int_{-\infty} ^{+\infty}\xi^4 \mu(d\xi ,d\tau ) ~~
\; ~~ \ \lambda_{20} = \int_{-\infty} ^{+\infty}\xi^2
\mu(d\xi,d\tau),
$$
where $\mu$ is the spectral measure of the stationary random field $
W(x,t)$. The density in (\ref{f:twin}) satisfies
$$
 p_{\mathbf{Y}(x,t)}(0) = (\lambda_{20})^{-1/2}\varphi\big(
 kx(\lambda_{20})^{-1/2}\big)(\lambda_{40})^{-1/2}\varphi\big(
 k(\lambda_{40})^{-1/2}\big).
$$
On the other hand
$$
\mathbf{Y}'(x,t) = \left( \begin{array}{cc}
  W_{xx}(x,t)  -k & W_{xt}(x,t)\\
  W_{xxx}(x,t) & W_{xxt}(x,t)
\end{array} \right).
$$
Under the condition $ \mathbf{Y}(x,t) =0$, one has
$$
|\det (\mathbf{Y}'(x,t))| = |W_{xt}(x,t)W_{xxx}(x,t)|.
$$
Computing the regression it turns out that the conditional
distribution of the pair $(W_{xt}(x,t),W_{xxx}(x,t))$ under the same
condition, is the one of two independent centered gaussian random
variables, with the following parameters:
\begin{align}
\mbox{expectation } \frac{\lambda_{31}}{\lambda_{40}} k &\mbox{ and
variance } \lambda_{22} -
\frac{\lambda^2_{31}}{\lambda_{40}},~ \mbox{for the first coordinate} \\
\mbox{expectation } \frac{\lambda_{40}}{\lambda_{20}} kx &\mbox{ and
variance } \lambda_{60} - \frac{\lambda^2_{40}}{\lambda_{20}},~
\mbox{for the second coordinate}
\end{align}
It follows that:
$$
 \E \Big(|\det (\mathbf{Y}'(x,t))| \s
\mathbf{Y}(x,t) =0\Big)  = G \Big(
\frac{\lambda_{31}}{\lambda_{40}} k, \sqrt{ \lambda_{22} -
\frac{\lambda^2_{31}}{\lambda_{40}}} \Big) .G \Big(
\frac{\lambda_{40}}{\lambda_{20}} kx, \sqrt{ \lambda_{60} -
\frac{\lambda^2_{40}}{\lambda_{20}}} \Big)
$$
Summing up:
\begin{multline}
 \frac{1}{T}\E\big( \mathcal{TW}(\R,T)\big)=
 \\
\frac{ 1}{\sqrt{\lambda_{40}}}  \varphi \left(\frac{ k}{\sqrt{\lambda_{40}}}\right)
G \Big(
\frac{\lambda_{31}}{\lambda_{40}} k, \sqrt{ \lambda_{22} -
\frac{\lambda^2_{31}}{\lambda_{40}}} \Big)
\int_\R G \Big(
\frac{\lambda_{40}}{\lambda_{20}} kx, \sqrt{ \lambda_{60} -
\frac{\lambda^2_{40}}{\lambda_{20}}} \Big)\frac{ 1}{\sqrt{\lambda_{20}}}
  \varphi \left(\frac{ kx}{\sqrt{\lambda_{20}}}\right)
  \\=
  \sqrt{ \frac{2}{\pi}}  \varphi \left(\frac{ k}{\sqrt{\lambda_{40}}}\right)
G \Big(\frac{\lambda_{31}}{\lambda_{40}} k, \sqrt{ \lambda_{22} -
\frac{\lambda^2_{31}}{\lambda_{40}}} \Big)
\frac{1}{k} \frac{\sqrt{\widetilde{ \lambda}_6 \lambda_{20} }
+\lambda_{40}}{\sqrt{ \lambda_{20}\lambda_{40}}}
\frac{\sqrt{\widetilde{ \lambda}_6}}{\widetilde{ \lambda}_6 +\lambda^2_{40} }
\end{multline}
  setting $\widetilde{ \lambda}_6 :=\lambda_{60} -
\frac{\lambda^2_{40}}{\lambda_{20}}$
 This result is equivalent to  formula (4.7)
  of \cite{lo1:lo1} (part III page 853).

\subsection{Specular points in two dimensions} \label{s:spec2d}
We consider at fixed time a random surface depending on two space
variables $x$ and $y$. The source of light is placed at $(0,0,h_1)$
and the observer is at $(0,0, h_2)$. The point $(x,y)$ is a specular
point if the normal vector $ n(x,y) = ( -W_x,-W_y,1)$ to the surface
at $(x,y)$ satisfies the following two conditions:
\begin{itemize}
  \item  the angles with the incident ray  $I = (
  -x,-y,h_1-W)$ and the reflected ray $R = (
  -x,-y,h_2-W)$ are equal (for
  short  the argument $(x,y)$ has been removed),
  \item it belongs to the plane generated by $I$ and $R$.
\end{itemize}
Setting $ \alpha_i =h_i -W$ and $ r_i = \sqrt{x^2 +y^2 +\alpha_i}$,
$i=1,2$, as in the one-parameter case we have:
 \begin{align}
W_x& = \frac{x}{x^2 +y^2}\frac{ \alpha _2 r_1 - \alpha_1 r_2}{r_2-r_1}, \notag
\\
W_y& = \frac{y}{x^2 +y^2}\frac{ \alpha _2 r_1 - \alpha_1 r_2}{r_2-r_1}.\label{spec:2:1}
\end{align}
When $h_1$ and $h_2$  are large, the system above can be approximated by
 \begin{align}
W_x& = kx \notag
\\
W_y& = ky,\label{spec:2:2}
\end{align}
under the same conditions as in dimension 1.\\

Next, we compute the expectation of $ SP_2(Q)$, the number of
approximate specular points in the sense of (\ref{spec:2:2}) that
are in a domain $Q$. In the remaining of this paragraph we limit our
attention to this approximation and to the case in which  $\{
W(x,y):(x,y)\in \R^2 \}$ is a centered
Gaussian stationary random field.\\

Let us define:
\begin{equation}\label{f:specd}
\mathbf{Y}(x,y) := \left( \begin{array}{c}
  W_x(x,y) -kx \\
W_{y} (x,y)-ky
\end{array} \right).
\end{equation}

Under very general conditions, for example on the spectral measure
of $\{W(x,y):x,y\in \R  \}$ the random field $\{Y(x,y):x,y\in \R \}$
satisfies the conditions of Theorem \ref{Rice2}, and we can write:
\begin{equation}\label{f:spec22}
\E\big( SP_2(Q)\big)  = \int_Q \E\big( |\det\mathbf{Y}'(x,y)|\big)
p_{\mathbf{Y}(x,y)}(\mathbf{0})\ dx dy,
\end{equation}
since for fixed $(x,y)$ the random matrix $Y'(x,y)$ and the random
vector $Y(x,y)$ are independent, so that the condition in the
conditional expectation can be erased.\\

The density in the
 \rhs of (\ref{f:spec22}) has the expression
\begin{equation}\aligned\label{f:d0}
p_{\mathbf{Y}(x,y)}(\mathbf{0})  &= p_{ (W_x,W_y)} (kx,ky)\\
&=
\frac{1}{2\pi}\frac{1}{\sqrt{\lambda_{20}\lambda_{02}-\lambda_{11}^2}
} \exp \Big[
-\frac{k^2}{2(\lambda_{20}\lambda_{02}-\lambda_{11}^2)}\big(\lambda_{02}
x^2-2\lambda_{11}xy+\lambda_{20}y^2 \big)\Big].
\endaligned
\end{equation}
To compute the expectation of the absolute value of the determinant
in the \rhs of (\ref{f:spec22}), which does not depend on $x,y$, we
use the method of \cite{be1:be1}. Set $ \Delta :=
\det\mathbf{Y}'(x,y)=( W_{xx} -k) (W_{yy} -k) -
W^2_{xy}$.\\

We have
\begin{equation} \label{f:d1}
\E(|\Delta|) = \E \left[\frac{2}{\pi} \int_0 ^{+\infty} \frac{1-
\cos( \Delta t)}{t^2} dt \right].
\end{equation}
Define
$$
h(t) :=  \E \left[ \exp \big( it [ ( W_{xx} -k) ( W_{yy} -k) -W_{xy}^2] \big)\right].
$$
Then
\begin{equation}\label{f:d2}
\E(|\Delta|) =\frac{2}{\pi}  \Big( \int_0 ^{+\infty} \frac{1-
\mathfrak{Re} [h(t)]}{t^2} dt \Big).
\end{equation}
To compute $h(t)$  we define
$$
A = \left( \begin{array}{ccc}
  0 & 1/2& 0\\
  1/2 & 0 & 0\\
  0& 0 & -1
\end{array}
\right)
$$
and  $ \Sigma$ the variance matrix  of $ W_{xx}, W_{yy},W_{x,y}$
$$
\Sigma : = \left( \begin{array}{ccc}
  \lambda_{40}  & \lambda_{22}& \lambda_{31}\\
  \lambda_{22} & \lambda_{04} & \lambda_{13}\\
  \lambda_{31}& \lambda_{13} & \lambda_{22}
\end{array}
\right).
$$
Let$\Sigma^{1/2}  A\Sigma^{1/2} = P \
diag(\Delta_1,\Delta_2,\Delta_3) P^T$  where $P$ is orthogonal. Then
by a diagonalization argument
\begin{multline} \label{m:1}
h(t) =  e^{itk^2}
\\
 \E\Big( \exp\big[ it \big(  ( \Delta_1 Z^2_1 -k (s_{11}  + s_{21}) Z_1 )
 + ( \Delta_2 Z^2_2 -k (s_{12}  + s_{22}) Z_2 )
 +( \Delta_3 Z^2_3 -k (s_{13}  + s_{23}) Z_3 )\big) \big] \Big)
\end{multline}
 where $ (Z_1,Z_2,Z_3)$ is standard normal  and $ s_{ij} $ are the entries of
 $\Sigma^{1/2} P^T$.

One can check that  if  $ \xi$ is  a standard normal variable and
$\tau, \mu$ are real constants, $\tau >0$:

 $$
  \E\big( e^{ i\tau ( \xi + \mu)^2 } \big) =(1-2i\tau)^{-1/2}e^{\frac{i\tau\mu^2}{(1-2i\tau)}}
  = \frac{1}{(1+4\tau ^2)^{1/4}}
  \exp\Big[\frac{-2\tau}{1+4\tau^2} +i\big(\varphi+\frac{\tau \mu ^2}{1+4\tau ^2}
  \big)\Big],
 $$
 where
 $$
 \varphi=\frac{1}{2}\arctan (2\tau ),~~0<\varphi <\pi/4.
 $$
Replacing in (\ref{m:1}), we obtain for $\mathfrak{Re} [h(t)]$ the
formula:

\begin{equation}\label{Reh(t)}
\mathfrak{Re} [h(t)]=\Big[
\prod_{j=1}^3\frac{d_j(t,k)}{\sqrt{1+4\Delta _j^2t^2}}\Big]\cos
\Big( \sum_{j=1}^3\big( \varphi _j (t)+k^2t\psi _j (t)\big)\Big)
\end{equation}
where, for $j=1,2,3$:
\begin{itemize}
\item $\displaystyle
d_j(t,k)=\exp\Big[-\frac{k^2t^2}{2}~\frac{(s_{1j}+s_{2j})^2}{1+4\Delta
_j^2t^2} \Big], $
\item $\displaystyle
 \varphi _j (t)=\frac{1}{2}\arctan (2\Delta _j t),~~0<\varphi _j<\pi/4,
 $
\item $\displaystyle
\psi _j (t)=\frac{1}{3}-t^2\frac{(s_{1j}+s_{2j})^2\Delta
_j}{1+4\Delta _j^2t^2}. $
\end{itemize}

Introducing these expressions in (\ref{f:d2})  and using
(\ref{f:d0}) we obtain a new formula which has the form of a rather
complicated integral. However, it is well adapted to numerical evaluation.\\

On the other hand, this formula allows us to compute the equivalent
as $k\rightarrow 0$ of the expectation of the total number of
specular points under the longuet-Higgins approximation. In fact, a
first order expansion of the terms in the integrand gives a somewhat
more accurate result, that we state as a theorem:
\begin{theo}\label{expspecdim2}
\begin{equation}\label{explonguetdim2}
\E \big(SP_2(\R^2)\big) = \frac{m_2}{k^2}+O(1)
\end{equation}
where
\begin{equation}\aligned \label{2dimlonghigg}
m_2&=\int_0^{+\infty}\frac{1-\big[\prod_{j=1}^3(1+4\Delta_j^2t^2)\big]^{-1/2}
\cos\big(\sum_{j=1}^3\varphi_j(t) \big)}{t^2}dt\\
&=\int_0^{+\infty}\frac{1-2^{-3/2} \big[\prod_{j=1}^3
\big(A_j\sqrt{1+A_j} \big)\big]\big( 1-B_1B_2-B_2B_3-B_3B_1\big)
}{t^2}dt,
\endaligned
\end{equation}
where
$$ A_j=A_j(t)=\big(1+4\Delta _j^2t^2
\big)^{-1/2},~~B_j=B_j(t)=\sqrt{(1-A_j)/(1+A_j)}.
$$
\end{theo}

Notice that $m_2$ only depends on the eigenvalues
$\Delta_1,\Delta_2,\Delta_3$ and is easily computed numerically.\\

In Flores and Le\'on \cite{flo:flo} a different approach was
followed in search of a formula for the expectation of the number of
specular points in the two-dimensional case, but their result is
only suitable for Montecarlo approximation.

\medskip

We now consider the \textbf{variance} of the total number of
specular points in two dimensions, looking for analogous results to
the one-dimensional case (i.e. Theorem \ref{asymptoticvar} and its
Corollary \ref{corcoefvar}), in view of their interest for
statistical applications. It turns out that the computations become
much more involved. The statements on variance and  speed of
convergence to zero of the coefficient of variation that we give
below include only the order of the asymptotic behavior in the
longuet-Higgins approximation, but not the constant. However, we
still consider them to be useful. If one refines the computations
one can give rough bounds on the generic constants in Theorem
\ref{varspecdim2} and Corollary \ref{coefvarspecdim2} on the basis
of additional
hypotheses on the random field.\\

We assume that the real-valued, centered, Gaussian stationary random
field $\{W(\mathbf{x}):\mathbf{x}\in \R^2\}$ has paths of class
$C^3$, the distribution of $W'(\mathbf{0})$ does not degenerate
(that is $\Var(W'(\mathbf{0}))$ is invertible). Moreover, let us
consider $W''(\mathbf{0})$, expressed in the reference system $xOy$
of $\R^2 $ as the $2\times 2$ symmetric centered Gaussian random
matrix:
$$W''(\mathbf{0})=\left(
    \begin{array}{cc}
      W_{xx}(\mathbf{0}) & W_{xy}(\mathbf{0}) \\
      W_{xy}(\mathbf{0}) & W_{yy}(\mathbf{0}) \\
    \end{array}
  \right)
$$
The function
$$
\mathbf{z}\rightsquigarrow \Delta (\mathbf{z})=\det \big[ \Var \big(
W''(0)\mathbf{z} \big) \big],
$$
defined on $\mathbf{z}=(z_1,z_2)^T \in \R^2$, is a non-negative
homogeneous polynomial of degree $4$ in the pair $z_1,z_2$. We will
assume the non-degeneracy condition:
\begin{equation}\label{nondegdelta}
\min \{\Delta (\mathbf{z}): \|\mathbf{z}\|=1 \}=\underline{\Delta
}>0.
\end{equation}

\begin{theo}\label{varspecdim2}
Let us assume that $\{W(\mathbf{x}):\mathbf{x}\in \R^2 \}$ satisfies
the above condtions and that it is also $\delta$-dependent, $\delta
>0 $, that is, $\E \big(W(\mathbf{x})W(\mathbf{y})\big)=0$
whenever $\|\mathbf{x}-\mathbf{y}\|>\delta.$\\

Then, for $k$ small enough:
$$
\Var \big(SP_2(\R^2) \big)~\leq ~ \frac{L}{k^2},
$$
where $L$ is a positive constant depending upon the law of the
random field.
\end{theo}

A direct consequence of Theorems \ref{expspecdim2} and
\ref{varspecdim2} is the following:

\begin{cor}\label{coefvarspecdim2}
Under the same hypotheses of Theorem \ref{varspecdim2}, for $k$
small enough, one has:
$$
\frac{\sqrt{\Var \big(SP_2(\R^2) \big)}}{\E\big(SP_2(\R^2)
\big)}~\leq ~ L_1k
$$
where $L_1$ is a new positive constant.
\end{cor}

\emph{Proof of Theorem \ref{varspecdim2}}.
For short, let us denote $T=SP_2( \R^2)$. We have:
\begin{equation}\label{vardim2}
\Var(T)=\E(T(T-1))+\E(T)-[\E(T)]^2
\end{equation}

We have already computed the equivalents as $k\rightarrow 0$ of the
second and third term in the right-hand side of (\ref{vardim2}). Our
task in what follows is to consider the first term.

The proof is performed along the same lines as the one of Theorem
\ref{asymptoticvar}, but instead of applying Rice formula for the
second factorial moment of the number of crossings of a
one-parameter random process, we need Theorem \ref{kricefield} for
dimension $d=2$. We write the factorial moment of order $m=2$ in the
form:

\begin{equation}\aligned
&\E(T(T-1))\\
&=\int\int_{\R^2 \times
\R^2}\E\Big(|\det\mathbf{Y}'(\mathbf{x})||\det\mathbf{Y}'(\mathbf{y})|\s~
\mathbf{Y}(\mathbf{x})=\mathbf{0},\mathbf{Y}(\mathbf{y})=\mathbf{0}
\Big) \notag
p_{\mathbf{Y}(\mathbf{x}),\mathbf{Y}(\mathbf{y})}(\mathbf{0},\mathbf{0})~d\mathbf{x}d\mathbf{y}\\
&=\int\int_{\|\mathbf{x}-\mathbf{y}\|>\delta}...~d\mathbf{x}d\mathbf{y}+\int\int_{\|\mathbf{x}-\mathbf{y}\|\leq
\delta}...~d\mathbf{x}d\mathbf{y}=J_1+J_2.\label{2dim2}
\endaligned
\end{equation}

For $J_1$ we proceed as in the proof of Theorem \ref{asymptoticvar},
using the $\delta$-dependence and the evaluations leading to the
statement of Theorem \ref{expspecdim2}. We obtain:
\begin{equation}\label{j1}
J_1=\frac{m_2^2}{k^4}+\frac{O(1)}{k^2}.
\end{equation}
Let us show that for small $k$,
\begin{equation}\label{j2}
J_2=\frac{O(1)}{k^2}.
\end{equation}
In view of (\ref{vardim2}), (\ref{explonguetdim2}) and
(\ref{j1}) this suffices to prove the theorem.\\

We do not perform all detailed computations. The key point consists
in evaluating the behavior of the integrand that appears in $J_2$
near the diagonal $\mathbf{x}=\mathbf{y}$, where the density
$p_{\mathbf{Y}(\mathbf{x}),\mathbf{Y}(\mathbf{y})}(\mathbf{0},\mathbf{0})$
degenerates and the conditional expectation tends to zero.\\

For the density, using the invariance under translations of the law
of $W'(\mathbf{x}):\mathbf{x}\in \R^2$, we have:
\begin{equation*}\aligned
p_{\mathbf{Y}(\mathbf{x}),\mathbf{Y}(\mathbf{y})}(\mathbf{0},\mathbf{0})&=
p_{W'(\mathbf{x}),W'(\mathbf{y})}(k\mathbf{x},k\mathbf{y})\\
&=
p_{W'(\mathbf{0}),W'(\mathbf{y-x})}(k\mathbf{x},k\mathbf{y})\\
&=p_{W'(\mathbf{0}),[W'(\mathbf{y-x})-W'(\mathbf{0})]}(k\mathbf{x},k
(\mathbf{y}-\mathbf{x})).
\endaligned
\end{equation*}
Perform the Taylor expansion, for small
$\mathbf{z}=\mathbf{y}-\mathbf{x} \in \R^2$:
$$
W'(\mathbf{z})=W'(\mathbf{0})+
W''(\mathbf{0})\mathbf{z}+O(\|\mathbf{z}\|^2).
$$
Using the non-degeneracy assumption (\ref{nondegdelta}) and the fact
that $W'(\mathbf{0})$ and $W''(\mathbf{0})$ are independent, we can
show that for $\mathbf{x},\mathbf{z} \in \R^2, \|\mathbf{z}\|\leq
\delta$:
$$
p_{\mathbf{Y}(\mathbf{x}),\mathbf{Y}(\mathbf{y})}(\mathbf{0},\mathbf{0})
\leq~\frac{C_1}{\|\mathbf{z}\|^2}\exp\big[-C_2 k^2
(\|\mathbf{x}\|-C_3 )^2 \big]
$$
where $C_1,C_2, C_3$ are positive constants.\\

Let us consider the conditional expectation. For each pair
$\mathbf{x},\mathbf{y}$ of different points in $\R^2$, denote by
$\mathbf{\tau }$ the unit vector
$(\mathbf{y}-\mathbf{x})/\|\mathbf{y}-\mathbf{x}\|$  and
$\mathbf{n}$ a unit vector orthogonal to $\mathbf{\tau }$. We denote
respectively by $\partial _{\mathbf{\tau }}\mathbf{Y},
\partial _{\mathbf{\tau }\mathbf{\tau }}\mathbf{Y},\partial _{\mathbf{n }}\mathbf{Y}$
the first and second partial derivatives of the random field in the
directions
given by $\mathbf{\tau }$ and $\mathbf{n } $.\\

Under the condition
$$
\mathbf{Y}(\mathbf{x})=\mathbf{0},\mathbf{Y}(\mathbf{y})=\mathbf{0}
$$
we have the following simple bound on the determinant, based upon
its definition and Rolle's Theorem applied to the segment
$[\mathbf{x},\mathbf{y}]=\{\lambda \mathbf{x}+(1-\lambda
)\mathbf{y}\}$:

\begin{equation}
\big|\det \mathbf{Y'}(\mathbf{x}) \big| \leq \|\partial
_{\mathbf{\tau }}\mathbf{Y}(\mathbf{x})\|\|\partial _{\mathbf{n
}}\mathbf{Y}(\mathbf{x})\|\leq
\|\mathbf{y}-\mathbf{x}\|\sup_{\mathbf{s}\in
[\mathbf{x},\mathbf{y}]}\|
\partial _{\mathbf{\tau }\mathbf{\tau }}\mathbf{Y}(\mathbf{s})\|\|\partial _{\mathbf{n
}}\mathbf{Y}(\mathbf{x})\|
\end{equation}
So,
\begin{equation*}\aligned
\E
&\Big(|\det\mathbf{Y}'(\mathbf{x})||\det\mathbf{Y}'(\mathbf{y})|\s~
\mathbf{Y}(\mathbf{x})=\mathbf{0},\mathbf{Y}(\mathbf{y})=\mathbf{0}
\Big)\\
&\leq \|\mathbf{y}-\mathbf{x}\|^2 \E \Big[\sup_{\mathbf{s}\in
[\mathbf{x},\mathbf{y}]}\|
\partial _{\mathbf{\tau }\mathbf{\tau }}\mathbf{Y}(\mathbf{s})\|^2\|\partial _{\mathbf{n
}}\mathbf{Y}(\mathbf{x})\|\|\partial _{\mathbf{n
}}\mathbf{Y}(\mathbf{y})\|
\Big|W'(\mathbf{x})=k\mathbf{x},W'(\mathbf{y})=k\mathbf{y} \Big]\\
&=\|\mathbf{z}\|^2 \E \Big[\sup_{\mathbf{s}\in
[\mathbf{0},\mathbf{z}]}\|
\partial _{\mathbf{\tau }\mathbf{\tau }}\mathbf{Y}(\mathbf{s})\|^2\|\partial _{\mathbf{n
}}\mathbf{Y}(\mathbf{0})\|\|\partial _{\mathbf{n
}}\mathbf{Y}(\mathbf{z})\|
\Big|W'(\mathbf{0})=k\mathbf{x},\frac{W'(\mathbf{z})-W'(\mathbf{0})}{\|\mathbf{z}\|}=k
\mathbf{\tau} \Big],
\endaligned
\end{equation*}
where the last equality is again a consequence of the stationarity
of the random field $\{ W(\mathbf{x}):\mathbf{x}\in \R^2 \}$.\\

At this point, we perform a Gaussian regression on the condition.
For the condition, use again Taylor expansion, the non-degeneracy
hypothesis and the independence of $W'(\mathbf{0})$ and
$W''(\mathbf{0})$. Then, use the finiteness of the moments of the
supremum of bounded Gaussian processes (see for example \cite{aw},
Ch. 2), take into account that $\|z\|\leq \delta $ to get the
inequality:
\begin{equation}\label{cotaespcond}
\E
\Big(|\det\mathbf{Y}'(\mathbf{x})||\det\mathbf{Y}'(\mathbf{y})|\s~
\mathbf{Y}(\mathbf{x})=\mathbf{0},\mathbf{Y}(\mathbf{y})=\mathbf{0}
\Big) \leq C_4~\|\mathbf{z}\|^2\big(1+k\|\mathbf{x}\| \big)^4
\end{equation}
where $C_4 $ is a positive constant. Summing up, we have the
following bound for $J_2$:

\begin{equation}\aligned
J_2 &\leq C_1C_4~ \pi \delta ^2 \int_{\R^2}\big(1+k\|\mathbf{x}\|
\big)^4\exp\big[-C_2 k^2 (\|\mathbf{x}\|-C_3 )^2 \big] d
\mathbf{x}\\
&= C_1C_4~ 2 \pi^2 \delta ^2 \int_{0}^{+\infty}\big(1+k\rho
\big)^4\exp\big[-C_2 k^2 (\rho -C_3 )^2 \big] \rho d\rho
\endaligned
\end{equation}
Performing the change of variables $w=k\rho$, (\ref{j2}) follows.
\findem

\section{The distribution of the normal to the level curve}  \label{s:normal}

Let us consider a modeling of the sea $ W(x,y,t)$ as a function
of two space variables and one time variable. Usual models are centered
Gaussian stationary with a particular form of the spectral measure $\mu$
that we discuss briefly below. We denote the covariance by $ \Gamma (x,y,t)=\E (W(0,0,0)W(x,y,t))$.\\

In practice, one is frequently confronted with the following situation:
several pictures of the sea on time over an interval $[0,T]$ are
stocked and some properties or magnitudes are observed. If the time
$T$ and the number of pictures are large, and if the process is
ergodic in time, the frequency of pictures that satisfy a certain
property will converge to the probability of this property to happen
at a fixed time.

Let us illustrate this with the angle of the normal to the level
curve at a point ``chosen at random''. We consider first the number
of crossings of a level $u$  by the process $W(\cdot,y,t)$ for fixed
$t$ and  $y$, defined as
$$
N_{[0,M_1]}^{W(\cdot, y,t)}(u)=\#\{ x :0\le x \le
M_1;\,W(x,y,t)=u\}.
$$
We are interested in computing the total number of crossings per
unit time when integrating over $ y \in [0,M_2]$   i.e.
\begin{equation}\label{hori}
\frac1T \int_0^T dt \int_0^{M_2}N_{[0,M_1]}^{W(\cdot,y,t)}(u)~dy.
\end{equation}
If the ergodicity assumption in time holds true, we can conclude that a.s.:
$$
\frac1T
\int_0^T dt \int_0^{M_2}N_{[0,M_1]}^{W(\cdot,y,t))}(u)~dy \to M_1
\E\big(N_{[0,M_1]}^{W(\cdot,0,0))}(u)\big)=\frac{M_1M_2}{\pi}\sqrt{\frac{\lambda_{200}}{\lambda_{000}}}e^{-\frac12\frac{u^2}{\lambda_{000}}},
$$
where
$$
\lambda_{abc}  =  \iint_{\R^3} \lambda_x^a\lambda_y^b\lambda_t^c
d\mu(\lambda_x,\lambda_y,\lambda_t)
$$
are the spectral moments.\\

Hence, on the basis of the quantity (\ref{hori}) for large $T$, one
can make inference about the value of certain parameters of the law
of the random field. In this example these are the spectral moments
$\lambda_{200}$ and $\lambda_{000}$.\\

If two-dimensional level information is available, one can work
differently because there exists an interesting relationship with Rice formula for level curves that we
explain in what follows.\\

We can write ($\textbf{x}=(x,y)$):
$$
W'(\mathbf{x},t)=||W'(\mathbf{x},t)||(\cos \Theta(\mathbf{x},t),\sin
\Theta(\mathbf{x},t))^T.
$$
 Instead of using  Theorem \ref{rice1}, we can use Theorem  \ref{rice4}, to
 write
\begin{align}
\E\big[
\int_0^{M_2}N_{[0,M_1]}^{W(\cdot,y,0)}(u)~dy\big]&=\E\big[\int_{{\cal
C}_Q(0,u)}|\cos\Theta(\mathbf{x},0)|~d\sigma_1\big] \notag
\\
&=\frac{\sigma_2(Q)}{\pi}\sqrt{\frac{\lambda_{200}}{\lambda_{000}}}\,e^{-\frac{u^2}
{2\lambda_{000}}},\label{f:chichi}
\end{align}
where $ Q = [0,M_1]\times[0,M_2] $. We have  a similar formula when
we consider  sections of the set $ [0,M_1]\times[0,M_2] $ in the
other direction. In fact (\ref{f:chichi}) can be generalized to
obtain the Palm distribution of the angle $ \Theta$.\\

Set $ h_{\theta_1,\theta_2}=\UN_{[
\theta_1\,,\,\theta_2]}$, and for
$-\pi\le\theta_1<\theta_2\le \pi$  define
\begin{align}
F(\theta_2)-F(\theta_1):&= \E\Big(\sigma_1(\{\mathbf{x}\in Q:\,
W(\mathbf{x},0)=u\, ;\, \theta_1\le \Theta(\mathbf{x},s)\le
\theta_2\}) \notag
\\
&= \E\Big(\int_{\mathcal{C}_Q(u,s)}h_{\theta_1,\theta_2}(
\Theta(\mathbf{x},s))d\sigma_1(\mathbf{x})ds\Big)\label{f:hh}
\\
&=\sigma_2(Q)\E[h_{\theta_1,\theta_2}(\frac{\partial_{y}W}
{\partial_{x}W})((\partial_{x}W)^2+(\partial_{y}W)^2)^{1/2}]\,\frac{\exp(-\frac{u^2}{2\lambda_{00}})}{\sqrt{
2\pi \lambda_{000}}}.\notag
\end{align}
Denoting $\Delta=\lambda_{200}\lambda_{020}-\lambda_{110}$  and
assuming $\sigma_2(Q)=1$ for ease of notation, we readily obtain
\begin{multline*}
 F(\theta_2)-F(\theta_1)
\\
=\frac{\,e^{-\frac{u^2}{2\lambda_{000}}}}
{(2\pi)^{3/2}(\Delta)^{1/2}\sqrt{\lambda_{000}}}\int_{\R^2}h_{\theta_1,\theta_2}(\Theta)
\sqrt{x^2+y^2}e^{-\frac1{2\Delta}(\lambda_{02}x^2-2\lambda_{11}xy+\lambda_{20}y^2)}dxdy
\\
=\frac{\,e^{-\frac{u^2}{2\lambda_{00}}}}{(2\pi)^{3/2}(\lambda_{+}\lambda_{-})^{1/2}\sqrt{\lambda_{000}}} \\
\int_0^\infty\int_{\theta_1}^{\theta_2}\rho^2
\exp(-\frac{\rho^2}{2\lambda_{+}\lambda_{-}}(\lambda_{+}\cos^2(\varphi-\kappa)+\lambda_{-}
\sin^2(\varphi-\kappa)))d\rho d\varphi
\end{multline*}
where $\lambda_{-}\le\lambda_{+}$ are the eigenvalues of the
covariance matrix of the random vector
$(\partial_{x}W(0,0,0),\partial_{y}W(0,0,0))$ and   $\kappa$ is
the angle of the eigenvector associated to $ \gamma^+$. Remarking
that the exponent  in the integrand can be written as\\  $
1/\lambda_-(1 -\gamma^2 sin^2( \varphi -\kappa)) $ with $\gamma^2 :=
1- \lambda_+/\lambda_-$  and that
$$
\int_0^{+\infty} \rho^2 \exp\Big( -\frac{H \rho^2}{2}\Big) =
\sqrt{\frac{\pi}{2H}}
$$
it is easy to get that
$$
F(\theta_2)-F(\theta_1) = (const) \int_{\theta_1}^{\theta_2} \big(1
-\gamma^2 sin^2( \varphi -\kappa)\big)^{-1/2} d \varphi.
$$
From this relation we get the density $g(\varphi)$ of  the Palm
distribution, simply by dividing by the total mass:
\begin{equation}\label{f:palm}
g(\varphi) = \frac{\big(1 -\gamma^2 \sin^2( \varphi
-\kappa)\big)^{-1/2}}{ \int_{-\pi}^{\pi} \big(1 -\gamma^2 \sin^2(
\varphi -\kappa)\big)^{-1/2}d \varphi.} =\frac{\big(1 -\gamma^2
\sin^2( \varphi -\kappa)\big)^{-1/2}}{ 4 \mathcal{K}(\gamma^2)} ,
\end{equation}
Here $ \mathcal{K}$ is the complete elliptic integral of the  first
kind. This density characterizes the distribution of the angle of
the normal at a point chosen ``at random'' on the level curve.\\

In the case of a random field which is isotropic  in $(x,y)$, we have
$\lambda_{200}=\lambda_{020}$ and moreover $\lambda_{110}=0$, so
that $g$ turns out to be the uniform density over the circle
(Longuet-Higgins says that over the contour the ``distribution''
of the angle is uniform (cf. \cite{lo2:lo2},  pp. 348)).\\


Let now  $\mathcal{W}=\{W(\mathbf{x},t):t \in \R^+,
\mathbf{x}=(x,y)\in \R^2 \}$ be a stationary zero mean Gaussian
random field modeling the height of the sea waves. It has the
following spectral representation:
$$W(x,y,t)=\int_{\Lambda}e^{i(\lambda_1x+\lambda_2y+\omega
t)}\sqrt{f(\lambda_1,\lambda_2,\omega)}dM(\lambda_1,\lambda_2,\omega),$$
where $\Lambda$ is the manifold
$\{\lambda_1^2+\lambda_2^2=\omega^4\}$ (assuming that the
acceleration of gravity $g$ is equal to 1) and $M$ is a random
Gaussian orthogonal measure defined on $\Lambda$ (see
\cite{kreso:kreso}). This leads to the following representation for
the covariance function
\begin{align*}
\Gamma(x,y,t) &=\int_{\Lambda}e^{i(\lambda_1x+\lambda_2y+\omega
t)}f(\lambda_1,\lambda_2,\omega)\sigma_2(dV)
\\
&=\int_{-\infty}^{\infty}\int_0^{2\pi}e^{i(\omega^2x\cos\varphi+\omega^2y\sin\varphi+\omega
t)}G(\varphi,\omega)d\varphi d\omega,
\end{align*}
where, in the second equation,
we made the change of variable $\lambda_1=\omega^2\cos\varphi$,
$\lambda_2=\omega^2\sin\varphi$ and
$G(\varphi,\omega)=f(\omega^2\cos\varphi,\omega^2\sin\varphi,\omega)2\omega^3$. The function $G$ is called the ``directional spectral function''. If
$G$ does not depend of $\varphi$ the random
field $W$ is isotropic in $x,y$.\\

Let us turn to ergodicity. For a given subset $Q$ of $\R^2$ and each $t$, let us define
$$\mathcal{A}_t=\sigma \{W(x,y,t):\tau >t \,;(x,y)\in
Q\}$$  and consider the $\sigma$-algebra of $t$-invariant events $ \mathcal{A}=\bigcap \mathcal{A}_t $. We assume that for each pair $(x,y)$, $\Gamma (x,y,t)\rightarrow 0$ as $t\rightarrow +\infty.$ It is well-known that under this condition, the $\sigma$-algebra $ \mathcal{A} $ is trivial, that is, it only contains events having probability zero or one (see for example \cite{cr:le}, Ch. 7).\\

This has the following important consequence in our context. Assume further that the set $Q$ has a smooth boundary and for simplicity, unit Lebesgue measure. Let us consider
$$Z(t)=\int_{\mathcal{C}_Q(u,t)}H\big(\mathbf{x},t\big)d\sigma_1(\mathbf{x}),
$$
where $ H\big(\mathbf{x},t\big)=\mathcal{H}\big( W(\mathbf{x},t),\nabla W (\mathbf{x},t )\big)$, where $\nabla W =(W_x,W_y)$ denotes gradient in the space variables and $\mathcal{H} $ is some measurable function such that the integral is well-defined. This is exactly our case in
(\ref{f:hh}). The process $\{Z(t):t\in \R\}$ is strictly stationary, and in our case has a finite mean and is
Riemann-integrable. By the Birkhoff-Khintchine ergodic theorem
(\cite{cr:le} page 151), a.s. as $T\rightarrow +\infty $,
$$
\frac1T\int_0^T Z(s)ds\to \E_{\mathcal{B}}[Z(0)],
$$
where $ \mathcal{B}$ is the $\sigma$-algebra of $t$-invariant
events associated to the process $Z(t)$. Since for each $t$, $Z(t)$ is $\mathcal{A}_t$-measurable, it follows that $ \mathcal{B} \subset \mathcal{A}$, so that $\E_{\mathcal{B}}[Z(0)]=\E[Z(0)] $.
On the other hand, Rice's formula yields (take into account that stationarity of $\mathcal{W} $ implies that $W(\mathbf{0},0)$ and $\nabla W(\mathbf{0},0)$ are independent):
\begin{align*}
\E[Z(0)]=
\E[\mathcal{H}\big( u,\nabla W(\mathbf{0},0)\big)||\nabla W(\mathbf{0},0)||]p_{W(\mathbf{0},0)}(u).
\end{align*} \medskip

We consider now the CLT. Let us define
$$
\mathcal{Z}(t)=\frac1t\int_0^t \big[ Z(s) -\E ( Z(0)) \big] ds,
$$

In order to compute second moments, we use Rice formula for integrals over level sets (cf. Theorem
\ref{rice4}), applied to the vector-valued random field
$$
X(\mathbf{x}_1,\mathbf{x}_2,s_1,s_2)= (W(\mathbf{x}_1,s_1),W(\mathbf{x}_2,s_2))^T.
$$
The level set can be written as:
$$
\mathcal{C}_{Q^2}(u,u)= \{(\mathbf{x}_1,\mathbf{x}_2)\in Q \times Q:\, X(\mathbf{x}_1,\mathbf{x}_2,s_1,s_2)=(u,u)\}
\qquad \mbox{ for } 0\le s_1\le t,\, 0\le s_2\le t.
$$
So, we get
$$
\Var \mathcal{Z}(t)=\displaystyle\frac2t\int_0^t (1-\frac st)I(u,s)ds,
$$
where
\begin{align*}
 I(u,s)=\int_{Q^2} &\E \Big[H(\mathbf{x}_1, 0) H(\mathbf{x}_2, s) \|\nabla W(\mathbf{x}_1,0)\|\| \nabla W(\mathbf{x}_2,s)\| \Big|~W(\mathbf{x}_1,0)=u\,;W(\mathbf{x}_2,s)=u \Big]\\
&\times p_{W(\mathbf{x}_1,0),W(\mathbf{x}_2,s)}(u,u) d\mathbf{x}_1 d\mathbf{x}_2~-  \Big(\E[\mathcal{H}\big( u,\nabla W(\mathbf{0},0)\big)||\nabla W(\mathbf{0},0)||]p_{W(\mathbf{0},0)}(u) \Big)^2
\end{align*}
Assuming that the given random field is time-$\delta$-dependent, that is,\\
$\Gamma(x,y,t)=0\, ~\forall ~(x,y) $, whenever $t>\delta$, we readily
obtain
$$
t\Var \mathcal{\mathcal{Z}}(t)\to2\int_0^{\delta}
I(u,s)ds:=\sigma^2(u) \ \ \mbox{ as }t\to\infty.
$$
Using now a variant of the Hoeffding-Robbins Theorem \cite{hoe:ro} for sums of
$\delta$-dependent random variables, we get the CLT:

$$
\sqrt t \mathcal{Z}(t)\Rightarrow N(0,\sigma^2(u)).
$$

%
\section{ Numerical computations} \label{s:num}

\subsection*{Validity of the approximation for the number of specular points} In the particular case of stationary processes we have compared
the exact expectation given by (\ref{spec1}) with the
approximation (\ref{f:spec2}).

In full generality the result depends on $ h_1,h_2, \lambda_4$ and
$\lambda_2$. After scaling, we can assume for example that $
\lambda_2=1$.

The main result is that, when $h_1 \approx h_2$, the approximation
(\ref{f:spec2}) is very sharp.
For example
 with the value $(100,100,3)$ for $  (h_1,h_2, \lambda_4) $,
 the expectation of the total number of specular
 points over $ \R$  is $138.2$; using the approximation (\ref{f:spectot})
 the result with the exact formula is
 around $ 2 .10^{-2}$ larger but it is almost hidden by  the precision of the
 computation of the integral.

 If we consider the case $(90,110,3)$, the results are respectively $136.81$ and $137.7$.

\begin{center}
\begin{figure}[h!]
\includegraphics[width =8 cm]{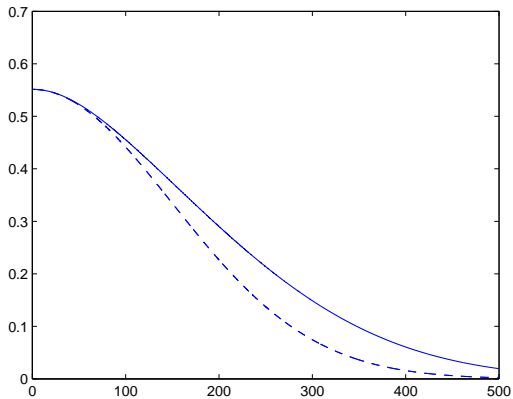}
\caption{ Intensity of specular points in the case $h_1= 100,h_2= 300,  \lambda_4=3$.
In solid line exact formula,  in dashed  line approximation (\ref{f:spec2}) }\label{fi:1}
\end{figure}
\end{center}

In the case $ (100,300,3)$,  the results differ significantly and
Figure \ref{fi:1}  displays the densities (\ref{spec1}) and
(\ref{f:spec2})

\subsection*{Effect of anisotropy on the distribution of the angle of the normal to the curve}
\begin{center}
\begin{figure}[h!]
\includegraphics[width =8 cm]{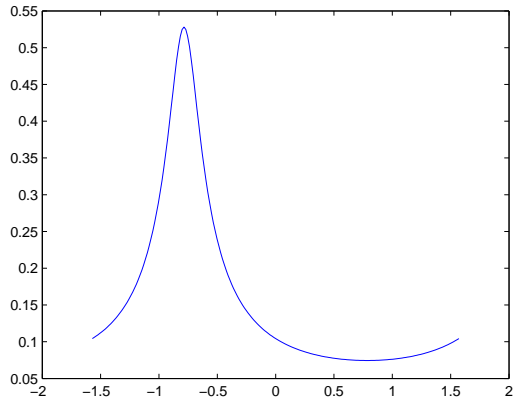}
\caption{ Density of the Palm  distribution of the angle of the  normal to the level curve in the case $\gamma =0.5$ and
$\kappa=\pi/4$ } \label{fi:2}
\end{figure}
\end{center}

We show the values of the density given by (\ref{f:palm}) in the
case of anisotropic processes $\gamma =0.5$ and $\kappa =\pi/4$.
Figure \ref{fi:2}  displays the densities  of the Palm distribution
of the angle showing  a large departure from the uniform
distribution.

\subsection*{Specular points in dimension 2}
\begin{center}
\begin{figure}[h!]
\includegraphics[width =8 cm]{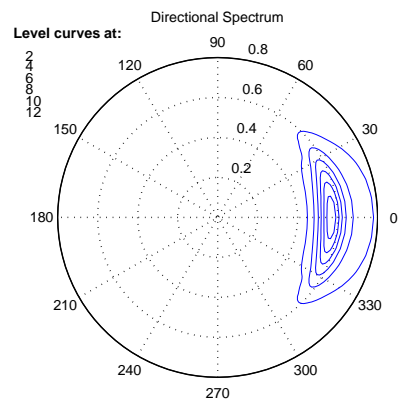}
\caption{ Directional  Jonswap spectrum  as obtained using the
default options of Wafo}\label{f:spectrum}
\end{figure}
\end{center}

We use a standard sea model   with a Jonswap spectrum  and spread
function $\cos(2 \theta)$. It corresponds to
 the default parameters of the Jonswap function of the toolbox WAFO \cite{wafo}.  The variance matrix of the gradient
 is equal to
 $$
  10^{ -4}\left(
   \begin{array}{cc}
     114 &0\\
      0& 81\\
   \end{array}
 \right)
 $$
 and the matrix $\Sigma$ of Section  \ref{s:spec2d}  is
 $$
 \Sigma =10^{ -4}  \left(
            \begin{array}{ccc}
              9& 3& 0  \\
              3& 11& 0 \\
              0& 0 & 3 \\
            \end{array}
          \right)
$$
The spectrum is presented in  Figure  \ref{f:spectrum}

\begin{center}
\begin{figure}[h!]
\includegraphics[width =8 cm]{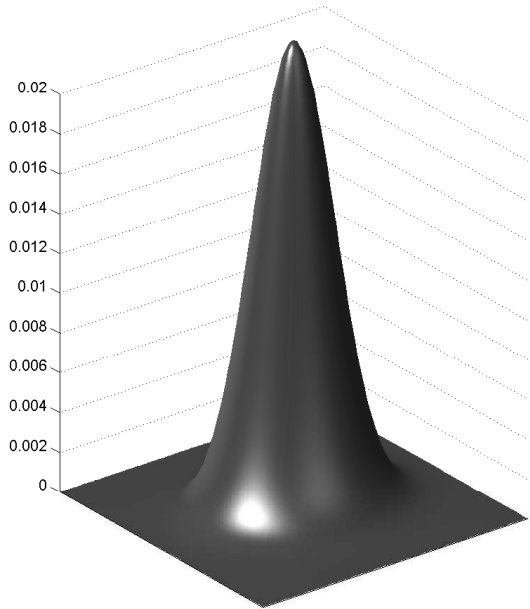}
\caption{ Intensity  function of the specular points  for  the
Jonswap spectrum}\label{f:spec2d}
\end{figure}
\end{center}

 The integrand  in (\ref{f:spec22})  is displayed in Figure \ref{f:spec2d} as a function of the two  space variables $x,y$. The value of the asymptotic parameter $m_2$ defining the expansion  on the expectation of the numbers of specular points, see(\ref{explonguetdim2}), is
 $2.527 10 ^{-3}$.

The Matlab programs  used for these computations are available at
\begin{verbatim}
http://www.math.univ-toulouse.fr/~azais/prog/programs.html
\end{verbatim}

\section{Application to dislocations of
wavefronts}\label{section4}

In this section we follow the article by Berry and Dennis
\cite{be1:be1}. As these authors, we are interested in dislocations
of wavefronts. These are lines in space or points in the plane where
the phase $\chi$, of the complex scalar wave
$\displaystyle\psi(\mathbf{x},t)=\rho(\mathbf{x},t)e^{i\chi(\mathbf{x},t)},$
is undefined, ($ \mathbf{x} = (x_1,x_2))$  is a two dimensional
space variable). With respect to light they are lines of darkness;
with respect to sound, threads of silence.\\

It will be convenient to express $\psi$ by means of its real and
imaginary parts:
$$
\psi(\mathbf{x},t)=\xi(\mathbf{x},t)+i\eta(\mathbf{x},t).
$$
Thus the
dislocations are the intersection of the two surfaces
$$
\xi(\mathbf{x},t)=0\qquad ~ \eta(\mathbf{x},t)=0.
$$
We assume an isotropic Gaussian model. This means that we will
consider the wavefront as an isotropic Gaussian field
$$
\psi({\bf x},t)=\int_{{\R}^2}\exp{(i[ \langle{\bf k}\cdot {\bf x} \rangle-c
|{\bf k}|t])}(\frac{\Pi(|{\bf k}|)}{|{\bf k}|})^{1/2}dW({\bf k}),
$$
 where,
${\bf k}=(k_1,k_2)$, $|{\bf k}|=\sqrt{k_1^2+k_2^2}$, $\Pi(k)$ is the
isotropic spectral density and $W=(W_1+iW_2)$ is a standard complex
orthogonal Gaussian measure on $\R^2$, with unit variance. Here
 we are interested only in $t=0$ and  we put
   $\xi({\bf x}):=\xi({\bf
x},0)$ and $\eta({\bf x}):=\eta({\bf x},0)$.\\

We have, setting $k = |{\bf k}|$
\begin{equation}\label{process1}
\xi({\bf x})=\int_{{\R}^2}\cos( \langle{\bf k}\cdot {\bf
x}\rangle)(\frac{\Pi(k)}{k})^{1/2}dW_1({\bf k})-\int_{{\R}^2} \sin(\langle{\bf
k}\cdot {\bf x}\rangle)(\frac{\Pi(k)}{k})^{1/2}dW_2({\bf k})
\end{equation}
and
\begin{equation}\label{process2}
\eta({\bf x})=\int_{{\R}^2}\cos(\langle {\bf k}\cdot {\bf
x}\rangle)(\frac{\Pi(k)}{k})^{1/2}dW_2({\bf k})+\int_{{\R}^2} \sin(\langle{\bf
k}\cdot {\bf x}\rangle)(\frac{\Pi(k)}{k})^{1/2}dW_1({\bf k})
\end{equation}
The covariances are
\begin{equation}\label{covariance}
\E\,[\xi({\bf x})\xi({\bf x}')]=\E\,[\eta({\bf x})\eta({\bf x}')]
= \rho(|\mathbf{x}-\mathbf{x}'|) :=
\int_0^\infty J_0(k|{\bf x}-{\bf x}'|)\Pi(k)dk
\end{equation}
where $J_{\nu}(x)$ is the Bessel function of the first kind of order
$\nu$. Moreover $\E\,[\xi({\bf r}_1)\eta({\bf r}_2)]=0.$

\subsubsection*{Three dimensional model}
 In the case of a three dimensional
Gaussian field, we have ${\bf
x}=(x_1,x_2,x_3)$, ${\bf
k}=(k_1,k_2,k_3)$,$k=|{\bf k}|=\sqrt{k_1^2+k_2^2+k_3^2}$ and
$$
\psi({\bf x})=\int_{{\R}^3}\exp\big(i[ \langle {\bf k}\cdot {\bf x} \rangle] \big)
(\frac{\Pi(k)}{k^2})^{1/2}dW({\bf k}).
$$
In this case, we write the covariances in the form:
\begin{eqnarray}\label{int1}
\E\,[\xi({\bf r}_1)\xi({\bf r}_2)]
=4\pi\int_0^\infty\frac{\sin(k|{\bf r}_1-{\bf
r}_2|)}{k|{\bf r}_1-{\bf r}_2|}\,\Pi(k)dk.
\end{eqnarray} The same
formula holds true for the process $\eta$ and also $ \E\,[\xi({\bf
r}_1)\eta({\bf r}_2)]=0$  for any $\bf {r}_1,\bf {r}_2$, showing
that the two coordinates are independent Gaussian fields .

\subsection{Mean length of dislocation curves,
mean number of dislocation points} {\bf Dimension 2:} Let us denote
$\{\mathbf{Z}(\mathbf{x}): \mathbf{x} \in \R^2\}$ a random field
with values in $\R^2$, with coordinates $\xi(\mathbf{x}),
\eta(\mathbf{x})$, which are two independent Gaussian stationary
isotropic random fields with the same distribution. We are
interested in the expectation of the number of dislocation points
 $$
 d_2 :=  \E [ \#\{\mathbf{x} \in S :\xi(\mathbf{x})
  =\eta(\mathbf{x})= 0\}],
 $$
where $S$ is a subset of the parameter space having area equal to
$1$.\\

Without loss of generality we may assume that $
\Var(\xi(\mathbf{x})) = \Var(\eta(\mathbf{x})) =1 $ and for the
derivatives  we set $ \lambda_2=  \Var (\eta_i(\mathbf{x}))=\Var
(\xi_i(\mathbf{x}))$, $ i=1,2$. Then, using stationarity and the
Rice formula (Theorem \ref{Rice2}) we get
$$
d_2 =\E[ |\det(\mathbf{Z}'(\mathbf{x})) |/ \mathbf{Z}(\mathbf{x})
=0] p_{\mathbf{Z}(\mathbf{x})}(0),
$$
The stationarity implies independence between
$\mathbf{Z}(\mathbf{x})$ and $ \mathbf{Z}'(\mathbf{x})$ so that the
conditional expectation above is in fact an ordinary expectation.
The entries of $\mathbf{Z}'(\mathbf{x})$ are four independent
centered Gaussian variables with variance $\lambda_2$, so that, up
to a factor, $|\det( \mathbf{Z}'(\mathbf{x})) |$ is the area of the
parallellogram generated by two independent standard Gaussian
variables in $\R^2$. One can easily show that the distribution of
this volume is the product of independent square roots of a $
\chi^2(2)$ and a $ \chi^2(1)$ distributed random variables. An
elementary calculation gives then: $\E[
|\det(\mathbf{Z}'(\mathbf{x})) |] = \lambda_2 $. Finally, we get
$$
d_2 = \frac{1}{2\pi } \lambda_2
$$
This quantity is equal to  $ \frac{K_2}{4\pi }$ in Berry and Dennis
\cite{be1:be1} notations, giving  their formula (4.6).\\

{\bf Dimension 3:} In the case, our aim is to compute
$$
d_3 =\E [ \mathcal{L}\{\mathbf{x} \in S :\xi(\mathbf{x})
=\eta(\mathbf{x})= 0\}]
$$
where $S$ is a subset of $\R^3 $ having volume equal to $1$ and $
\mathcal{L}$ is the length of the curve.  Note that $d_3$ is denoted
by  $d$ \cite{be1:be1}. We use the same notations and remarks except
that the form of the Rice's formula is (cf.  Theorem \ref{rice3})
$$
d_3  = \frac{1}{2\pi} \E[ (\det \mathbf{Z}'(\mathbf{\mathbf{x}})
\mathbf{Z}'(\mathbf{\mathbf{x}})^T)^{1/2} ].
$$
Again
$$
\E[ ( \det ( \mathbf{Z}'(\mathbf{x})
\mathbf{Z}'(\mathbf{x})^T)^{1/2}] = \lambda_2 \E( V),
$$
where $V$ is the surface area of the parallelogram generated by two
standard Gaussian variables in $\R^3$. A similar method to compute
the expectation of this random area gives:
 $$
 \E (V) = \E(\sqrt{\chi^2(3) }) \times \E(\sqrt{\chi^2(2 ) }) =  \frac{4}{\sqrt{2\pi}}
\sqrt{\frac{\pi}{2}} =2
$$
Leading eventually to
$$
d_3 =\frac{ \lambda_2}{\pi}.
$$
In Berry and Dennis' notations \cite{be1:be1} this last quantity is
denoted by $\frac{ k_2}{3\pi}$ giving their formula (4.5).

\subsection{Variance }

In this section we limit ourselves to dimension $\mathbf{2}$. Let
$S$ be again a measurable subset of $\R^2$ having Lebesgue measure
equal to $1$. The computation of the variance of the number of
dislocations points is performed using Theorem \ref{kricefield} to
express
$$
\E \big(N^\mathbf{Z}_S(\mathbf{0})\big(N^\mathbf{Z}_S(\mathbf{0})-1
\big) \big) = \int_{S^2}
A\mathbf{(s}_1,\mathbf{s}_2) d\mathbf{s}_1 d\mathbf{s}_2.
$$
We assume that $\{\mathbf{Z}(\mathbf{x}):\mathbf{x} \in \R^2\}$
satisfies the hypotheses of Theorem \ref{kricefield} for $m=2$.
Then use
$$
\Var\big( N^\mathbf{Z}_S(\mathbf{0})\big)=\E
\big(N^\mathbf{Z}_S(\mathbf{0})\big(N^\mathbf{Z}_S(\mathbf{0})-1
\big) \big)+d_2-d_2^2.
$$
Taking into account that the law of the random field is invariant
under translations and orthogonal transformations of $\R^2$, we have
$$
A(\mathbf{s}_1,\mathbf{s}_2)  =  A\big((0,0),(r,0) \big) = A(r)  \  \ \mbox{ whith } r= \| \mathbf{s}_1-\mathbf{s}_2\|,
$$

The  Rice's function $ A(r))$ ,  has two intuitive interpretations. First it can be viewed as
$$
A(r)= \lim_{\epsilon \to 0} \frac{1}{ \pi^2
\epsilon^4}\E\big[N\big( B( (0,0), \epsilon)\big)\times N\big( B(
(r,0), \epsilon)\big)\big].
$$
Second it is
 the density of the
Palm distribution (a generalization Horizontal window conditioning
of \cite{cr:le}) of the number of zeroes of $\mathbf{Z}$ per unit
of surface, locally around the point $(r,0)$ given that there is a
zero at
$(0,0)$.

\noindent
$A(r)/d_2^2$ is called ``correlation  function" in \cite{be1:be1}. \medskip

To compute $A(r)$, we put  $\xi_1,\xi_2,\eta_1\eta_2$
for the partial derivatives of $\xi, \eta$ with respect to first and
second coordinate.\\

and
\begin{align}
A(r)&= \E \big[ |\det  \mathbf{Z}'(0,0) \det \mathbf{Z}'(r,0)| \s
\mathbf{Z}(0,0) = \mathbf{Z}(r,0)=\mathbf{0}_2 \big] p_{\mathbf{
Z}(0,0), \mathbf{Z}(r,0)}(\mathbf{0}_4)\notag
\\
&=\E \big[ \big|\big(\xi_1\eta_2-\xi_2\eta_1\big)(0,0)\big(
\xi_1\eta_2-\xi_2\eta_1\big)(r,0)\big) \big| \Big|\mathbf{Z}(0,0) =
\mathbf{Z}(r,0)=\mathbf{0}_2\big]  \notag
\\
& \ \ \ \  \ \ \ \  \ \ \ \  \ \ \ \  p_{\mathbf{ Z}(0,0),\label{corre}
\mathbf{Z}(r,0)}(\mathbf{0}_4)
\end{align}
where $ \mathbf{0}_p$ denotes the null vector in dimension $p$.\\

The density is easy to compute
$$
p_{ \mathbf{Z}(0,0), \mathbf{Z}(r,0)}(\mathbf{0}_4) = \frac{1}{(2\pi)^2(1-\rho^2(r))},
\mbox{ where } \displaystyle\rho(r)=\int_0^\infty J_0(kr)\Pi(k)dk.
$$
We use now the same device as above to compute the conditional
expectation of the modulus of the product of determinants, that is
we write:
\begin{equation} \label{sauce:menthe}
|w| = \frac{1}{\pi} \int_{-\infty}  ^{ + \infty}  (1- \cos(wt)  t ^{-2} dt.
\end{equation}
and also the same notations as in \cite{be1:be1}

$$
\left\{
\begin{array}{c}
   C := \rho(r) \\
  E = \rho'(r)\\
  H= -E/r\\
  F =-\rho"(r) \\
 F_0=-\rho"(0)
\end{array} \right.
$$
The regression formulas imply that the conditional variance matrix
of the vector
$$
\mathbf{ W}=\Big( \xi_1(\mathbf{0}),\xi_1(r,0) ,
\xi_2(\mathbf{0}),\xi_2(r,0) ,\eta_1(\mathbf{0}),\eta_1(r,0),
\eta_2(\mathbf{0}),\eta_2(r,0) \Big),
$$
is given by
$$
\Sigma = Diag \Big[ \mathcal{A}, \mathcal{B},\mathcal{A}, \mathcal{B}\Big ]
$$
with
$$
\mathcal{A} = \left(\begin{array}{cc}
  F_0 -\frac{E^2}{1-C^2}& F -\frac{E^2 C}{1-C^2} \\
  F -\frac{E^2 C}{1-C^2}&F_0 -\frac{E^2}{1-C^2}
\end{array}  \right)
$$
$$
\mathcal{B} = \left(\begin{array}{cc}
  F_0 & H \\
  H&F_0
\end{array}  \right)
$$
Using formula (\ref{sauce:menthe}) the expectation we have to compute
 is equal to
 \begin{multline}\label{f:dbint}
 \frac{1}{\pi^2} \int_{_-\infty}  ^{ + \infty}  dt_1 \int_{- \infty}  ^{ + \infty} dt_2 t_1^{-2}
 t_2^{-2} \Big[
 1 - \frac 1 2  T(t_1,0)  - \frac 1 2  T(-t_1,0)-  \frac 1 2 T(0,t_2) -  \frac 1 2 T(0,-t_2)\\
 + \frac 1 4  T( t_1,t_2) + \frac 1 4  T( -t_1,t_2) + \frac 1 4  T( t_1,-t_2) + \frac 1 4  T( -t_1,-t_2) \Big]
\end{multline}
where
$$
T(t_1,t_2) = \E \big[ \exp\big(i ( w_1 t_1 +w_2t_2) \big) \Big]
$$
with
$$
w_1 =\xi_1(\mathbf{0})\eta_2(\mathbf{0}) -\eta_1(\mathbf{0})\xi_2(\mathbf{0}) = \mathbf{W}_1 \mathbf{W}_7 -\mathbf{W}_3 \mathbf{W}_5
$$
$$
w_2 =\xi_1(r,0)\eta_2(r,0) -\eta_1(r,0)\xi_2(r,0) = \mathbf{W}_2 \mathbf{W}_8 -\mathbf{W}_4 \mathbf{W}_6.
$$
$ T( t_1,t_2) = \E \big( \exp(i \mathbf{W}^T \mathcal{H }\mathbf{W}) \big)$  where $\mathbf{W}$
 has the distribution $N( 0, \Sigma)$  and
$$
 \mathcal{H }= \left[ \begin{array}{cccc}
   0 & 0& 0 & \mathcal{ D} \\
   0 & 0& -\mathcal{ D} & 0 \\
   0 & -\mathcal{ D}& 0 & 0 \\
  \mathcal{ D} & 0 & 0 & 0 \
 \end{array} \right  ],
 $$

$$
\mathcal{ D} = \frac{1}{2}\left[ \begin{array}{cc}
   t_1 & 0\\
   0 & t_2
\end{array} \right  ].
$$
A standard diagonalization argument shows that
$$
 T( t_1,t_2) = \E \big( \exp(i \mathbf{W}^T \mathcal{H }\mathbf{W}) \big)
 = \E\big( \exp( i \sum_{j=1} ^8 \lambda_j \xi_j^2)  \big),
$$
where the $ \xi_j$'s are independent with standard normal
 distribution and the $\lambda_j$ are the eigenvalues of
$ \Sigma ^{1/2} \mathcal{H} \Sigma ^{1/2}$. Using the characteristic
 function of the $ \chi^2(1)$
distribution:
\begin{equation}\label{chi}
\E \big( \exp(i \mathbf{W}^T \mathcal{H} \mathbf{W}) \big)  =  \prod_{j=1} ^8  ( 1-2i \lambda_j)^{ -1/2}.
\end{equation}
 Clearly
 $$
\Sigma ^{1/2}  = Diag \Big[ \mathcal{A}^{1/2}, \mathcal{B}^{1/2},\mathcal{A}^{1/2},
\mathcal{B}^{1/2}\Big ]
$$
and
$$
\Sigma ^{1/2} \mathcal{H}\Sigma ^{1/2} = \left[ \begin{array}{cccc}
   0 & 0& 0 & \mathcal{ M} \\
   0 & 0& -\mathcal{ M^T} & 0 \\
   0 & -\mathcal{ M}& 0 & 0 \\
  \mathcal{ M^T} & 0 & 0 & 0 \
 \end{array} \right  ]
 $$
with $ \mathcal{M }=\mathcal{A}^{1/2} \mathcal{D} \mathcal{B}^{1/2}$.

Let  $\lambda$  be an eigenvalue of $ \Sigma ^{1/2} \mathcal{H}
\Sigma ^{1/2}$ It is easy to check that $\lambda^2$ is an eigenvalue
of $\mathcal{M }\mathcal{M}^T$. Respectively if $ \lambda_1^2$  and
$\lambda^2_2$ are the  eigenvalues of $\mathcal{M }\mathcal{M}^T$,
those  of $\Sigma ^{1/2} \mathcal{H}\Sigma ^{1/2}$ are $  \pm
\lambda_1$(twice) and $  \pm \lambda_2$ (twice).

Note that $ \lambda^2_1$ and $\lambda^2_2$ are the eigenvalues of
 $\mathcal{M} \mathcal{M}^T =\mathcal{A}^{1/2}
\mathcal{D} \mathcal{B}\mathcal{D} \mathcal{A}^{1/2}$ or
equivalently, of  $\mathcal{D} \mathcal{B}\mathcal{D} \mathcal{A}$.
Using (\ref{chi})
$$
\E \big( \exp(i \mathbf{W}^T \mathcal{H} \mathbf{W}) \big)  = \big( 1 +4 (\lambda^2_1 + \lambda^2_2)
 +16 \lambda^2_1  \lambda^2_2\big)^{ -1}
 =\big( 1 +4 tr ( \mathcal{D} \mathcal{B}\mathcal{D} \mathcal{A})
   +16 \det  ( \mathcal{D} \mathcal{B}\mathcal{D} \mathcal{A})\big)^{ -1}
 $$
where
$$ \mathcal{D} \mathcal{B}\mathcal{D} \mathcal{A} =
\frac{1}{4}\left[ \begin{array}{cc}
     t_1^2 F_0  (F_0 -\frac{E^2}{1-C^2} ) + t_1 t_2 H (F -\frac{E^2 C}{1-C^2})
                     &  t_1^2 F_0 (F -\frac{E^2 C}{1-C^2})   + t_1 t_2 H  (F_0 -\frac{E^2}{1-C^2})
                     \\
   t_1 t_2 H (F_0 -\frac{E^2}{1-C^2})  + t_2^2F_0  (F -\frac{E^2 C}{1-C^2})
                 &  t_1 t_2 H (F -\frac{E^2 C}{1-C^2}) + t_2^2F_0(F_0 -\frac{E^2}{1-C^2})
\end{array} \right  ]
$$
So,
\begin{align}
 4 tr( \mathcal{D} \mathcal{B}\mathcal{D} \mathcal{A})
&= (t_1^2  + t_2^2)  F_0(F_0 -\frac{E^2}{1-C^2}) + 2t_1 t_2 H (F -\frac{E^2 C}{1-C^2})
\\
16 \det( \mathcal{D} \mathcal{B}\mathcal{D} \mathcal{A})
 &=t_1^2 t_2^2 \big[ F_0^2 -H^2\big] \big[(F_0 -\frac{E^2}{1-C^2})^2 -(F -\frac{E^2 C}{1-C^2})^2\big]
\end{align}
giving
\begin{multline}T( t_1,t_2) =
\E \big( \exp(i \mathbf{W}^T \mathcal{H} \mathbf{W}) \big)
\\
 =\Big( 1 +(t_1^2  +t_2^2)  F_0(F_0 -\frac{E^2}{1-C^2}) + 2t_1 t_2 H (F -\frac{E^2 C}{1-C^2})
\\
+
t_1^2 t_2^2 \big[ F_0^2 -H^2\big] \big[(F_0 -\frac{E^2}{1-C^2})^2 -(F -\frac{E^2 C}{1-C^2})^2\big]
\Big)^{ -1}
\end{multline}

Performing the change of variable $t' = \sqrt{A_1} t$ with $ A_1=F_0(F_0 -\frac{E^2}{1-C^2})$ the integral
(\ref{f:dbint}) becomes
\begin{multline}\label{int22}
\frac{A_1}{\pi^2} \int_{-\infty}  ^{ + \infty}  dt_1 \int_{- \infty}  ^{ + \infty} dt_2
 t_1^{-2}
 t_2^{-2}
 \\
  \bigg[ 1 - \frac{1}{ 1+t_1^2} \frac{1}{ 1+t_2^2} + - \frac12 \Big\{\frac{1}
 { 1+ (t_1^2+t_2^2) -2 A_2 t_1 t_2 +  t_1^2 t_2^2 Z } + \frac{1}
 { 1+ (t_1^2+t_2^2) +2 A_2 t_1 t_2 +  t_1^2 t_2^2 Z } \Big\} \bigg]
 \\
 = \frac{A_1}{\pi^2} \int_{-\infty}  ^{ + \infty}  dt_1 \int_{- \infty}  ^{ + \infty} dt_2
 t_1^{-2}
 t_2^{-2}
 \\
  \bigg[ 1 - \frac{1}{ 1+t_1^2} - \frac{1}{ 1+t_2^2}  + \frac{ 1+(t_1^2+t_2^2)  +  t_1^2 t_2^2 Z  }
 {\Big( 1+ (t_1^2+t_2^2) +  t_1^2 t_2^2 Z \Big)^2-4 A_2^2  t_1^2 t_2^2}  \bigg]
 \end{multline}
 where
$$
\left\{
\begin{array}{l}
   A_2 =        \frac{ H}{F_0} \frac{ F( 1-C^2) - E^2C}{F_0( 1-C^2) -E^2 }
   \\
  Z = \frac{F_0^2-H^2}{F_0^2} \Big[ 1-(F -\frac{E^2 C}{1-C^2})^2. (F_0 -\frac{E^2}{1-C^2})^{-2}
  \Big].
\end{array} \right.
$$
In this form, and up to a sign change, this result is equivalent to Formula (4.43) of \cite{be1:be1} (note that $A_2^2 = Y$ in \cite{be1:be1}).

In order to compute the integral (\ref{int22}), first we obtain
$$\int_{-\infty}^{\infty}\frac1{t_2^2}\big[1-\frac1{1+t_2^2}\big]dt_2=\pi.
$$
We split the other term into two integrals, thus we have for the first one\\
\\
$\displaystyle\frac12\int_{-\infty}^{\infty}\frac1{t_2^2}\big[\frac1{1+(t_1^2+t_2^2)-2A_2t_1t_2+t_1^2t_2^2Z}
-\frac1{1+t_1^2}\big]dt_2$
\begin{align*}&=-\frac1{2(1+t_1^2)}\int_{-\infty}^{\infty}\frac1{t_2^2}\frac{(1+t_1^2Z)
t_2^2-2A_2t_1t_2}{1+t_1^2-2A_2t_1t_2+(1+t_1^2Z)t_2^2}dt_2
\\
&=-\frac1{2(1+t_1^2)}\int_{-\infty}^{\infty}\frac1{t_2^2}\frac{t_2^2-2Z_1t_1t_2}{t_2^2-2Z_1t_1t_2+Z_2}dt_2=I_1,
\end{align*}
where $Z_2=\frac{1+t_1^2}{1+Zt_1^2}$ and $Z_1=\frac{A_2}{1+Zt_1^2}$.\\
Similarly for the second integral we get\\
\\
$\displaystyle\frac12\int_{-\infty}^{\infty}\frac1{t_2^2}\big[\frac1{1+(t_1^2+t_2^2)+2A_2t_1t_2+t_1^2t_2^2Z}
-\frac1{1+t_1^2}\big]dt_2$
$$
=-\frac1{2(1+t_1^2)}\int_{-\infty}^{\infty}\frac1{t_2^2}\frac{t_2^2+2Z_1t_1t_2}{t_2^2+2Z_1t_1t_2+Z_2}\,dt_2=I_2
$$
\begin{align*}I_1+I_2&=-\frac1{2(1+t_1^2)}\int_{-\infty}^{\infty}\frac1{t_2^2}
\big[\frac{t_2^2-2Z_1t_1t_2}{t_2^2-2Z_1t_1t_2+Z_2}+
\frac{t_2^2+2Z_1t_1t_2}{t_2^2+2Z_1t_1t_2+Z_2}\big]dt_2\\
&=-\frac1{(1+t_1^2)}\int_{-\infty}^{\infty}\frac{t_2^2+(Z_2-4Z_1^2t_1^2)}{t_2^4+2(Z_2-2Z^2_1t_1^2)t_2^2+Z_2^2}\,dt_2\\
&=-\frac1{(1+t_1^2)}
\frac{\pi(Z_2-2Z_1^2t_1^2)}{Z_2\sqrt{(Z_2-Z_1^2t_1^2)}} .
\end{align*}
In the third line we have used the formula provided by the method of
residues. In fact, if the polynomial $X^2-SX+P$ with $P>0$ has not
root in $[0,\infty)$, then
$$
\int_{-\infty}^{\infty}\frac{t^2-\gamma}{t^4-St^2+P}\,dt=\frac{\pi}{\sqrt{P(-S+2\sqrt
P)}}(\sqrt P-\gamma).
$$ In our case $\gamma=-(Z_2-4Z_1^2t_1^2)$,
$S=-2(Z_2-2Z^2_1t_1^2)$ and $P=Z_2^2$. \\
Therefore we get

$$
A(r)=  \frac{A_1}{4 \pi^3(1-C^2)}  \int_{-\infty}^{\infty} \frac1{t_1^2}\big[1-\frac1{(1+t_1^2)}
\frac{(Z_2-2Z_1^2t_1^2)}{ Z_2\sqrt{(Z_2-Z_1^2t_1^2)}}\big]dt_1.
$$

\subsection*{Acknowledgement}  This work has received financial support   from European Marie Curie Network SEAMOCS.

{}

\end{document}